\documentclass{article}
\usepackage{geometry}

\pdfoutput=1

\usepackage{amsmath,amssymb,amsthm}

\usepackage[utf8]{inputenc}
\usepackage[T1]{fontenc}
\usepackage{ae}
\usepackage{aecompl}

\usepackage{caption}
\usepackage{subcaption}
\usepackage{tikz}

\usepackage[breaklinks]{hyperref}

\theoremstyle{plain}
\newtheorem{theorem}{Theorem}
\newtheorem{lemma}[theorem]{Lemma}
\newtheorem{corollary}[theorem]{Corollary}
\newtheorem{proposition}[theorem]{Proposition}

\newtheorem{observation}[theorem]{Observation}

\theoremstyle{definition}
\newtheorem{definition}[theorem]{Definition}
\newtheorem{example}[theorem]{Example}

\theoremstyle{remark}
\newtheorem{remark}[theorem]{Remark}

\newcommand{\R}{\mathbb{R}}
\newcommand{\N}{\mathbb{N}}
\newcommand{\C}{\mathbb{C}}

\newcommand{\bigO}{\mathcal{O}}
\newcommand{\smallO}{o}
\newcommand{\asyOp}{\mathcal{A}}
\newcommand*\fring[3]{\R[[#1]]^{#2}_{#3}}
\newcommand*\cring[3]{\C[[#1]]^{#2}_{#3}}
\newcommand*\G[3]{\Gamma^{#2}_{#3}\left(#1\right)}

\author{Michael Borinsky}

\title{Generating Asymptotics \\for factorially divergent sequences}

\date{}

\begin{document}
\maketitle
\begin{abstract}
The algebraic properties of formal power series, whose coefficients show factorial growth and admit a certain well-behaved asymptotic expansion, are discussed. It is shown that these series form a subring of $\R[[x]]$. This subring is also closed under composition and inversion of power series. An `asymptotic derivation' is defined which maps a power series to the asymptotic expansion of its coefficients. Product and chain rules for this derivation are deduced. With these rules asymptotic expansions of the coefficients of implicitly defined power series can be obtained. The full asymptotic expansions of the number of connected chord diagrams and the number of simple permutations are given as examples. 
\end{abstract}
\tableofcontents
\section{Introduction}
This article%
\footnote{%
An extended abstract of this article appeared as a contribution to FPSAC 2017 \cite{borinsky2016generating}.%
} 
is concerned with real sequences $f_n$, which admit an asymptotic expansion for large $n$ of the form,
\begin{align} \label{eqn:first_example} f_n = \alpha^{n+\beta} \Gamma(n+\beta) \left( c_0 + \frac{c_1}{\alpha(n+\beta-1)} + \frac{c_2}{\alpha^2(n+\beta-1)(n+\beta-2)} + \cdots \right), \end{align}
for some $\alpha \in \R_{>0}$, $\beta \in \R$ and $c_k \in \R$.
Sequences of this type appear in many enumeration problems, which deal with coefficients of factorial growth. For instance, certain subclasses of permutations and graphs of fixed valence show this behaviour \cite{albert2003enumeration,bender1978asymptotic}. 
Furthermore, there are countless examples where \textit{perturbative expansions} of physical quantities admit asymptotic expansions of this kind \cite{bender1969anharmonic,le2012large,dunne2012resurgence}.

The restriction to this specific class of power series is inspired by the work of Bender. In \cite{bender1975asymptotic} he analyzed the asymptotic behaviour of the coefficients of the composition of a power series, which has mildly growing coefficients, with a power series, which has rapidly growing coefficients. Here, Bender's results are extended into a complete algebraic framework. This is achieved by making heavy use of generating functions in the spirit of the `analytic combinatorics' or `symbolic method' approach \cite{flajolet2009analytic,bergeron1998combinatorial,wilf2013generatingfunctionology}. The key step is to interpret the \textit{coefficients of the asymptotic expansion as another power series}.

The resulting framework bears many resemblances to the theory of resurgence, which was established by Jean Ecalle \cite{ecalle1981fonctions}. Resurgence assigns a special role to power series whose coefficients grow factorially, as they offer themselves to be Borel transformed. For instance, it can be used to assign a unique function to such a factorially divergent sequence. This function could be interpreted as the sequence' generating function. Moreover, resurgence provides a promising approach to cope with divergent perturbative expansions in physics. Its application to these problems is an active field of research \cite{alvarez2004langer,dunne2012resurgence,aniceto2011resurgence}.

During a conversation with David Sauzin it became plausible that the presented methods can also be derived from resurgence. In fact, the formalism can be seen as a toy model of resurgence's \textit{calcul diff\'erentiel \'etranger} \cite[Vol. 1]{ecalle1981fonctions} also called \textit{alien calculus} \cite[II.6]{mitschi2016divergent}. This toy model is unable to fully reconstruct functions from asymptotic expansions, but does not rely on analytic properties of Borel transformed functions and therefore offers itself for combinatorial applications. A detailed and illuminating account on resurgence theory is given in Sauzin's review \cite[Part II]{mitschi2016divergent}. 

\subsection{Statement of results} 
Power series whose coefficients have a well-behaved asymptotic expansion, as in eq.\ \eqref{eqn:first_example}, form a subring of $\R[[x]]$, which will be denoted as $\fring{x}{\alpha}{\beta}$. This subring is also closed under composition and inversion of power series. A linear map, $\asyOp^\alpha_\beta:\fring{x}{\alpha}{\beta}\rightarrow \R[[x]]$, can be defined which \textit{maps a power series to the asymptotic expansion of its coefficients}. A natural way to define such a map is to associate the power series $\sum_{n=0}^\infty c_n x^n$ to the series $\sum_{n=0}^\infty f_n x^n$ related as in eq.\ \eqref{eqn:first_example}. This map turns out to be a \textit{derivation} that means it fulfills a \textit{product rule}
\begin{align*} &\text{with $f,g\in \fring{x}{\alpha}{\beta}$}& (\asyOp^\alpha_\beta (f \cdot g))(x) &= f(x) (\asyOp^\alpha_\beta g)(x) + g(x) (\asyOp^\alpha_\beta f)(x) \\ &\text{and a \textit{chain rule}, }& (\asyOp^\alpha_\beta (f \circ g))(x) &= f'(g(x)) (\asyOp^\alpha_\beta g)(x) + \left(\frac{x}{g(x)} \right)^{\beta} e^{\frac{\frac{1}{x} - \frac{1}{g(x)}}{\alpha}} (\asyOp^\alpha_\beta f ) (g(x)), \end{align*}
where $(f \cdot g)(x) = f(x) g(x)$ and $(f \circ g)(x) = f(g(x))$. In the second line it is required that $g_0=0$ and $g_1=1$.
These statements will be derived from elementary properties of the $\Gamma$ function.

Note that the chain rule involves a peculiar correction term if the coefficients $f_n$ of the power series $f(x)$ have a non-trivial asymptotic expansion. It is obvious that the chain rule cannot be as simple as the ordinary chain rule for differentiation. For general $f,g\in \fring{x}{\alpha}{\beta}$: $(\asyOp^\alpha_\beta (f \circ g))(x) \neq f'(g(x)) (\asyOp^\alpha_\beta g)(x)$. Otherwise, the reasonable requirement that the coefficients of the generating function $g(x)=x$ have a trivial asymptotic expansion, $(\asyOp^\alpha_\beta g)(x) = 0$, would imply that all $f\in \fring{x}{\alpha}{\beta}$ have trivial asymptotic expansions.  

In Sections~\ref{sec:prereq}-\ref{sec:dgls} the derivation ring $\fring{x}{\alpha}{\beta}$ will be described in detail and the main Theorem~\ref{thm:chainrule}, which establishes the chain rule for the asymptotic derivation, will be proven. 
The formalism can be applied to calculate the asymptotic expansions of the coefficients of implicitly defined power series. This procedure is similar to the extraction of the derivative of an implicitly defined function using the implicit function theorem. We will use it in Section~\ref{sec:applications} to give the full asymptotic expansions of the number of \textit{connected chord diagrams} and the full asymptotic expansions of the number of \textit{simple permutations}.

\subsection{Notation}
A (formal) power series $f \in \R[[x]]$ will be denoted in the usual `functional' notation $f(x) = \sum_{n=0}^\infty f_n x^n$. The coefficients of a power series $f$ will be expressed by the same symbol with the index attached as a subscript $f_n$ or with the coefficient extraction operator $[x^n]f(x) = f_n$. Ordinary (formal) derivatives are expressed as $f'(x) = \sum_{n=0}^\infty n f_n x^{n-1}$. 
The (Cauchy) product of two power series $f,g$ will be expressed either as $f \cdot g$, $(f \cdot g)(x)$ or $f(x) g(x)$ depending on the context. Correspondingly, we will switch freely between the different notations $f \circ g$, $(f \circ g)(x)$ and $f( g(x))$ for the composition of two power series. The ring of power series, restricted to expansions of functions which are analytic at the origin, or equivalently power series with non-vanishing radius of convergence, will be denoted as $\R\{x\}$. The $\bigO$-notation will be used: $\bigO(a_n)$ denotes the set of all sequences $b_n$ such that $\limsup _{n\rightarrow \infty} | \frac{ b_n }{a_n} | < \infty$ and $\smallO(a_n)$ denotes all sequences $b_n$ such that $\lim _{n\rightarrow \infty} \frac{ b_n }{a_n} = 0$. Equations of the form $a_n = b_n + \bigO(c_n)$ are to be interpreted as statements $a_n-b_n \in \bigO(c_n)$ as usual. See \cite{bender1974asymptotic} for an introduction to this notation. Tuples of non-negative integers will be denoted by bold letters $\textbf{t} = (t_1,\ldots,t_L) \in \N_0^L$. The notation $|\textbf{t}|$ will be used as a short form for $\sum_{l=1}^L t_l$. We will consider the binomial coefficient ${ a \choose n}$ to be defined for all $a\in \R$ and $n\in \N_0$ by ${ a \choose n} := [x^n] (1+x)^a$.

The only non-standard notation that will be used to improve the readability of lengthy expressions is the abbreviation $\G{n}{\alpha}{\beta} := \alpha^{n+\beta} \Gamma(n+\beta)$.
\section{Prerequisites}
\label{sec:prereq}
We will start by defining the subset of power series whose coefficients have well-behaved asymptotic expansions:
\begin{definition}
\label{def:Fpowerseries}
For given $\alpha \in \R_{>0}$ and $\beta \in \R$ let $\fring{x}{{\alpha}}{\beta}$ be the subset of $\R[[x]]$, such that $f \in \fring{x}{{\alpha}}{\beta}$ if and only if there exists a sequence of real numbers $(c_{k}^f)_{k\in \N_0}$, which fulfills
\begin{align} \label{eqn:basic_asymp_C} f_n &= \sum _{k=0}^{R-1} c_{k}^f \G{n-k}{\alpha}{\beta} + \bigO\left(\G{n-R}{\alpha}{\beta}\right) && \forall R \in \N_0, \end{align}
where $\G{n}{\alpha}{\beta} = \alpha^{n+\beta} \Gamma(n+\beta)$.
\end{definition}
\begin{observation}
$\fring{x}{{\alpha}}{\beta}$ is a linear subspace of $\R[[x]]$.
\end{observation}
\begin{observation}
The sequence $(c_{k}^f)_{k\in \N_0}$ is unique for each fixed $f \in \fring{x}{{\alpha}}{\beta}$. The coefficients can be calculated iteratively using the explicit formula 
$c^f_K = \lim \limits_{n \rightarrow \infty} \frac{f_n - \sum _{k=0}^{K-1} c_{k}^f \G{n-k}{\alpha}{\beta} }{ \G{n-K}{\alpha}{\beta}}$ for all $K\in \N_0$.
\end{observation}
Both these properties follow immediately from Definition~\ref{def:Fpowerseries}.
\begin{remark}
The expression in eq.\ \eqref{eqn:basic_asymp_C} represents an asymptotic expansion or Poincar\'e expansion with the asymptotic scale $\alpha^{n+\beta}\Gamma(n+\beta)$ \cite[Ch. 1.5]{bruijn1970asymptotic}.
\end{remark}
\begin{remark}
The subspace $\fring{x}{\alpha}{\beta}$ includes all (real) power series whose coefficients only grow exponentially: $\R\{x\} \subset \fring{x}{\alpha}{\beta}$. These with all other series with coefficients, which are in $\smallO(\G{n-R}{\alpha}{\beta})$ for all $R\in \N_0$, have an asymptotic expansion of the form in eq.\ \eqref{eqn:basic_asymp_C} with all $c_{k}^f=0$.
\end{remark}
\begin{remark}
    Definition~\ref{def:Fpowerseries} implies that if $f\in \fring{x}{\alpha}{\beta}$, then 
\begin{align*} f_n \in \bigO\left(\G{n}{\alpha}{\beta}\right)=\bigO\left(\alpha^{n} \Gamma(n+\beta)\right). \end{align*}
Accordingly, the power series in $\fring{x}{\alpha}{\beta}$ are a subset of \textit{Gevrey-1} sequences \cite[Ch XI-2]{hsieh2012basic}. Being \textit{Gevrey-1} is not sufficient for a power series to be in $\fring{x}{\alpha}{\beta}$. For instance, a sequence which behaves for large $n$ as $f_n = n! (1+ \frac{1}{\sqrt{n}}+ \bigO(\frac{1}{n}))$ is \textit{Gevrey-1}, but not in $\fring{x}{\alpha}{\beta}$ for any pair $(\alpha, \beta)$.
\end{remark}
\begin{remark}
In resurgence theory further restrictions on the allowed power series are imposed, which ensure that the Borel transformations of the sequences have proper analytic continuations or are `endless continuable' \cite[II.6]{mitschi2016divergent}. These restrictions are analogous to the requirement that, apart from $f_n$, also the sequence $c_{k}^f$ has to have a well-behaved asymptotic expansion. The coefficients of this asymptotic expansion are also required to have a well-behaved asymptotic expansion and so on. These kinds of restrictions will not be necessary for the presented algebraic considerations, which are aimed at combinatorial applications.
\end{remark}

The central theme of this article is to \textit{interpret the coefficients $c_{k}^f$ of the asymptotic expansion as another power series}. 
In fact, Definition~\ref{def:Fpowerseries} immediately suggests to define the following map:
\begin{definition}
\label{def:basic_asymp_definition}
Let $\asyOp^\alpha_\beta: \fring{x}{\alpha}{\beta} \rightarrow \R[[x]]$ be the map that associates a power series $\asyOp^\alpha_\beta f \in \R[[x]]$ to every power series $f \in \fring{x}{\alpha}{\beta}$ such that
\begin{align} \label{eqn:basic_asymp} (\asyOp^\alpha_\beta f)(x) = \sum_{k=0}^\infty c_{k}^f x^k, \end{align}
with the coefficients $c_{k}^f$ from Definition~\ref{def:Fpowerseries}.
\end{definition}
\begin{observation}
$\asyOp^\alpha_\beta$ is linear.
\end{observation}
\begin{remark}
In Proposition~\ref{prop:derivation} it will be proven that $\asyOp^\alpha_\beta$ is a derivation. We will adopt the usual notation for derivations and consider $\asyOp^\alpha_\beta$ to act on everything to its right.
\end{remark}
\begin{remark}
In the realm of resurgence such an operator is called \textit{alien derivative} or \textit{alien operator} \cite[II.6]{mitschi2016divergent}.
\end{remark}
\begin{remark}
$\asyOp^\alpha_\beta$ is clearly not injective. For instance, $\R\{x\} \subset \ker \asyOp^\alpha_\beta$.
\end{remark}
\begin{example}
The power series $f\in \R[[x]]$ associated to the sequence $f_n = n!$ clearly fulfills the requirements of Definition~\ref{def:Fpowerseries} with $\alpha=1$ and $\beta=1$. Therefore, $f \in \fring{x}{1}{1}$ and $(\asyOp^1_1 f)(x) = 1$. 
\end{example}

The asymptotic expansion in eq.\ \eqref{eqn:basic_asymp_C} is normalized such that shifts in $k$, $c_k^f \rightarrow c_{k-m}^f$, can be absorbed by shifts in $\beta$, $\beta \rightarrow \beta +m$. More specifically,
\begin{proposition}
\label{prop:betashift}
For all $m \in \N_0$
\begin{align*} f\in \fring{x}{\alpha}{\beta}\text{ if and only if } f\in \fring{x}{\alpha}{\beta+m}\text{ and }\asyOp^\alpha_{\beta+m} f \in x^m \R[[x]]. \end{align*}
If either holds, then $x^m \left( \asyOp^\alpha_{\beta} f\right)(x) = \left( \asyOp^\alpha_{\beta+m} f\right)(x)$.
\end{proposition}

\begin{proof}
Because $\G{n}{\alpha}{\beta} = \alpha^{n-m+\beta+m} \Gamma(n-m+\beta+m) = \G{n-m}{\alpha}{\beta+m}$, the following two relations between $f_n$ and $c^f_k$ are equivalent,
\begin{align} \label{eqn:shifteq1} f_n &= \sum _{k=0}^{R-1} c_{k}^f\G{n-k}{\alpha}{\beta} + \bigO\left(\G{n-R}{\alpha}{\beta}\right) && \forall R\in \N_0 \\ \label{eqn:shifteq2} f_n &= \sum _{k=m}^{R'-1} c_{k-m}^f\G{n-k}{\alpha}{\beta+m} + \bigO\left(\G{n-R'}{\alpha}{\beta+m}\right) && \forall R'\geq m. \end{align}
Eq.\ \eqref{eqn:shifteq1} follows from $f\in \fring{x}{\alpha}{\beta}$ by Definition~\ref{def:Fpowerseries}. In that case, eq.\ \eqref{eqn:shifteq2} implies that $f\in \fring{x}{\alpha}{\beta+m}$ and that $\left( \asyOp^\alpha_{\beta+m} f\right)(x) = \sum_{k=m}^\infty c_{k-m}^f x^k= x^m\left( \asyOp^\alpha_{\beta} f\right)(x) \in x^m \R[[x]]$ by Definition~\ref{def:basic_asymp_definition}.

If $f\in \fring{x}{\alpha}{\beta+m}$ and $\asyOp^\alpha_{\beta+m} f \in x^m \R[[x]]$, then we can write the asymptotic expansion of $f$ in the form of eq.\ \eqref{eqn:shifteq2}. Eq.\ \eqref{eqn:shifteq1} implies $f\in \fring{x}{\alpha}{\beta}$.
\end{proof}

By analogous reasoning, we can absorb shifts in $n$, $f_n \rightarrow f_{n+m}$, in eq.\ \eqref{eqn:basic_asymp_C} by shifts in $\beta$, $\beta \rightarrow \beta +m$.
\begin{proposition}
\label{prop:betashiftlow}
For all $m \in \N_0$
\begin{align*} f\in \fring{x}{\alpha}{\beta} \cap x^m \R[[x]] \text{ if and only if } \frac{f(x)}{x^m} \in \fring{x}{\alpha}{\beta+m}. \end{align*}
If either holds, then $\left( \asyOp^\alpha_{\beta} f\right)(x) = \left( \asyOp^\alpha_{\beta+m} \frac{f(x)}{x^m} \right)(x)$.
\end{proposition}
\begin{proof}
Because $\G{n+m}{\alpha}{\beta} = \alpha^{n+m+\beta} \Gamma(n+m+\beta) = \G{n}{\alpha}{\beta+m}$, the following two relations between $f_n$ and $c_k^f$ are equivalent,
\begin{align} \label{eqn:shiftloweq1} f_n &= \sum _{k=0}^{R-1} c_{k}^f\G{n-k}{\alpha}{\beta} + \bigO\left(\G{n-R}{\alpha}{\beta}\right) && \forall R\in \N_0 \\ \label{eqn:shiftloweq2} f_{n+m} &= \sum _{k=0}^{R-1} c_{k}^f\G{n-k}{\alpha}{\beta+m} + \bigO\left(\G{n-R}{\alpha}{\beta+m}\right) && \forall R\in \N_0. \end{align}
Eq.\ \eqref{eqn:shiftloweq1} follows from $f\in \fring{x}{\alpha}{\beta}$. Because $f\in x^m \R[[x]]$, we have $\frac{f(x)}{x^m} = \sum_{n=0}^\infty f_{n+m} x^{n}\in \R[[x]]$. Eq.\ \eqref{eqn:shiftloweq2} then implies that $\frac{f(x)}{x^m} \in \fring{x}{\alpha}{\beta+m}$ and by Definition~\ref{def:basic_asymp_definition}, $\left( \asyOp^\alpha_{\beta} f\right)(x) = \left( \asyOp^\alpha_{\beta+m} \frac{f(x)}{x^m}\right)(x)$. 

If $\frac{f(x)}{x^m} \in \fring{x}{\alpha}{\beta+m} \subset \R[[x]]$, then $f\in x^m \R[[x]]$ and eq.\ \eqref{eqn:shiftloweq2} holds for the coefficients of $f$, which implies $f\in \fring{x}{\alpha}{\beta}$ by eq.\ \eqref{eqn:shiftloweq1} and Definition~\ref{def:Fpowerseries}.
\end{proof}

From Proposition~\ref{prop:betashift}, it follows that $\fring{x}{\alpha}{\beta} \subset \fring{x}{\alpha}{\beta+m}$ for all $m\in \N_0$. It will be convenient to only work in the spaces $\fring{x}{\alpha}{\beta}$ with $\beta > 0$ and to use Proposition~\ref{prop:betashift} to verify that the subspaces $\fring{x}{\alpha}{\beta-m}$ inherit all relevant properties from $\fring{x}{\alpha}{\beta}$.
The advantage is that, with $\beta > 0$, it is easier to express uniform bounds on the remainder terms in eq.\ \eqref{eqn:basic_asymp_C}. The following definition will provide a convenient notation for these bounds.
\begin{definition}
\label{defn:basic_bigC_estimate}
For $\alpha,\beta \in \R_{>0}$ and $R \in \N_0$, let $\rho^\alpha_{\beta,R}: \fring{x}{\alpha}{\beta} \rightarrow \R_+$ be the map
\begin{align} \label{eqn:basic_bigC_estimate} \rho^\alpha_{\beta,R}(f) = \max_{0\leq K\leq R} \sup_{n\geq K} \frac{\left| f_n - \sum _{k=0}^{K-1} c_{k}^f \G{n-k}{\alpha}{\beta} \right|}{\G{n-K}{\alpha}{\beta}}, \end{align}
with the coefficients $c_k^f$ as in Definition~\ref{def:Fpowerseries}.
\end{definition}
It follows directly from Definition~\ref{def:Fpowerseries} that the quantity $\rho^\alpha_{\beta,R}(f)$ is finite. 
Eq.\ \eqref{eqn:basic_bigC_estimate} can be translated into bounds for the coefficients $f_n$ and the $c_{k}^f$:
\begin{observation}
\label{obs:smallc_bigC_estimate}
If $\alpha,\beta \in \R_{>0}$ and $R\in \N_0$, then for all $f \in \fring{x}{{\alpha}}{\beta}$ and $n,K \in \N_0$ with $K \leq R$ as well as $n \geq K$, 
\begin{align} \left| f_n - \sum _{k=0}^{K-1} c_{k}^f \G{n-k}{\alpha}{\beta} \right| &\leq \rho^\alpha_{\beta,R}(f) \G{n-K}{\alpha}{\beta} &&\text{ and } &&|c^f_K| \leq \rho^\alpha_{\beta,R}(f). \end{align}
\end{observation}
\begin{remark}
It can be verified using linearity and the triangle inequality that the maps $\rho^\alpha_{\beta,R}$ form a family of norms on all spaces $\fring{x}{\alpha}{\beta}$ where $\beta>0$. Moreover, these norms will turn out to be submultiplicative up to equivalence (see Proposition~\ref{prop:submultiplicative}). However, we will not make direct use of any topological properties of the spaces $\fring{x}{{\alpha}}{\beta}$ in this article.
\end{remark}

\section{Elementary properties of sums over $\Gamma$ functions}
The following lemma forms the foundation for most of the conclusions that will follow. It provides an estimate for sums of $\Gamma$ functions. Moreover, it ensures that the subspace $\fring{x}{\alpha}{\beta}$ of formal power series corresponds to a subset of a large class of sequences studied by Bender \cite{bender1975asymptotic}.
From another perspective the lemma can be seen as an entry point to resurgence, which bypasses the necessity for analytic continuations and Borel transformations.

\begin{lemma}
\label{lmm:gammasum}
If $\alpha,\beta \in \R_{>0}$, then 
\begin{align} \label{eqn:centersumexpl2} \sum_{m=0}^{n} \G{m}{\alpha}{\beta} \G{n-m}{\alpha}{\beta} &\leq (2+\beta)\G{0}{\alpha}{\beta} \G{n}{\alpha}{\beta} && \forall n \in \N_0. \end{align}
\end{lemma}
\begin{proof}
Recall that $\G{n}{\alpha}{\beta} = \alpha^{n+\beta}\Gamma(n+\beta)$ and that $\Gamma:\R_{>0}\rightarrow \R_{>0}$ is a \textit{log-convex} function. If $\beta \in \R_{>0}$, then the functions $\Gamma(m+\beta)$ and $\Gamma(n-m+\beta)$ are also log-convex functions in $m$ on the interval $[0,n]$, as log-convexity is preserved under shifts and reflections. 
Furthermore, log-convexity is closed under multiplication. This implies that $\G{m}{\alpha}{\beta}\G{n-m}{\alpha}{\beta} = \alpha^{n+2\beta} \Gamma(m+\beta)\Gamma(n-m+\beta)$ is a log-convex function in $m$ on the interval $[1,n-1] \subset [0,n]$. A convex function always attains its maximum on the boundary of its domain. Accordingly,
$\G{m}{\alpha}{\beta}\G{n-m}{\alpha}{\beta} \leq \G{1}{\alpha}{\beta}\G{n-1}{\alpha}{\beta}$ for all $m \in [1,n-1]$.
This way, the sum $\sum_{m=0}^{n}\G{m}{\alpha}{\beta}\G{n-m}{\alpha}{\beta}$ can be estimated after stripping off the two boundary terms:
\begin{align} \label{eqn:sumGestimate} \sum_{m=0}^{n} \G{m}{\alpha}{\beta}\G{n-m}{\alpha}{\beta} &\leq 2 \G{0}{\alpha}{\beta}\G{n}{\alpha}{\beta} + (n-1) \G{1}{\alpha}{\beta}\G{n-1}{\alpha}{\beta} && \forall n \geq 1. \end{align}
It follows from $n \Gamma(n) = \Gamma(n+1)$ that
$ \G{1}{\alpha}{\beta}\G{n-1}{\alpha}{\beta} = \frac{\beta}{n-1+\beta} \G{0}{\alpha}{\beta}\G{n}{\alpha}{\beta} $ for all $n\geq1$.
Because $n-1+\beta \geq n-1$, substituting this into eq.\ \eqref{eqn:sumGestimate} implies the inequality in eq.\ \eqref{eqn:centersumexpl2} for all $n \geq 1$. The remaining case $n=0$ is trivially fulfilled.
\end{proof}

\begin{corollary}
\label{crll:centersum}
If $\alpha,\beta \in \R_{>0}$ and $R \in \N_0$ are kept fixed, then there exists a constant $C\in \R$ such that
\begin{align} \label{eqn:centersum} \sum_{m=R}^{n-R} \G{m}{\alpha}{\beta}\G{n-m}{\alpha}{\beta} &\leq C \G{n-R}{\alpha}{\beta} && \forall n \geq 2R. \end{align}
\end{corollary}
\begin{proof}
Recall that $\G{m+R}{\alpha}{\beta} = \G{m}{\alpha}{\beta+R}$.
We can shift the summation variable to rewrite the left hand side of eq.\ \eqref{eqn:centersum} as 
\begin{gather*} \sum_{m=0}^{n-2R} \G{m+R}{\alpha}{\beta}\G{n-m-R}{\alpha}{\beta} = \sum_{m=0}^{n-2R} \G{m}{\alpha}{\beta+R}\G{n-2R-m}{\alpha}{\beta+R} \\ \leq (2+\beta+R)\G{0}{\alpha}{\beta+R} \G{n-2R}{\alpha}{\beta+R}, \end{gather*}
where we applied Lemma~\ref{lmm:gammasum} with the substitutions $\beta \rightarrow \beta+R$ and $n\rightarrow n-2R$. Because  $\G{n-2R}{\alpha}{\beta+R} = \G{n-R}{\alpha}{\beta}$ the statement follows.
\end{proof}

\begin{corollary}
\label{crll:exponential_inO}
If $\alpha,\beta \in \R_{>0}$, $C \in \R$ and $P \in \R[m]$ is some polynomial in $m$, then 
\begin{align} \sum_{m=R}^{n} C^m P(m) \G{n-m}{\alpha}{\beta} \in \bigO(\G{n-R}{\alpha}{\beta}) && \forall R \in \N_0. \end{align}
\end{corollary}
\begin{proof}
There is a constant $C' \in \R$ such that $|C^m P(m)|$ is bounded by $C' \G{m}{\alpha}{\beta}$ for all $m\in \N_0$. 
Therefore, Corollary~\ref{crll:centersum} ensures that 
\begin{gather*} \sum_{m=R}^{n-R} C^m P(m) \G{n-m}{\alpha}{\beta} \leq C'\sum_{m=R}^{n-R} \G{m}{\alpha}{\beta} \G{n-m}{\alpha}{\beta} \in \bigO(\G{n-R}{\alpha}{\beta}). \end{gather*}
The remainder 
$\sum_{m=n-R+1}^{n} C^m P(m) \G{n-m}{\alpha}{\beta}= \sum_{m=0}^{R-1} C^{n-m} P(n-m) \G{m}{\alpha}{\beta}$ is obviously in $\bigO(\G{n-R}{\alpha}{\beta})$.
\end{proof}

\section{A derivation for asymptotics}

\begin{proposition}
    If $\alpha \in \R_{>0}$, $\beta \in \R$ and $f,g\in \fring{x}{\alpha}{\beta}$, then 
\begin{itemize}
\item
    The product $(f \cdot g)(x) = f(x) g(x)$ belongs to $\fring{x}{\alpha}{\beta}$.
\item
\label{prop:derivation}
$\asyOp^\alpha_\beta$ is a derivation, that means it respects the product rule
\begin{align} \label{eqn:leipniz} ( \asyOp^\alpha_\beta ( f \cdot g ))(x) &= f(x) (\asyOp^\alpha_\beta g)(x) + g(x) (\asyOp^\alpha_\beta f)(x). \end{align}
\end{itemize}
\end{proposition}

\begin{corollary}
\label{crll:chain_for_product}
If $g^{1}, \ldots, g^{L} \in \fring{x}{\alpha}{\beta}$, then $\prod_{l=1}^L g^{l}(x) \in \fring{x}{\alpha}{\beta}$ and 
\begin{gather} \left( \asyOp^\alpha_\beta \left( \prod_{l=1}^L g^{l}(x) \right)\right)(x) = \sum_{l=1}^L \left(\prod_{\substack{m=1 \\ m\neq l}}^L g^{m}(x)\right) ( \asyOp^\alpha_\beta g^{l})(x).   \end{gather}
\end{corollary}
\begin{proof}
Proof by induction in $L$ using the product rule.
\end{proof}
\begin{corollary}
\label{crll:chain_for_product_pow}
If $g^{1}, \ldots, g^{L} \in \fring{x}{\alpha}{\beta}$ and $\textbf{t} = (t_1, \ldots, t_L) \in \N_0^L$, then $\prod_{l=1}^L (g^{l}(x))^{t_l} \in \fring{x}{\alpha}{\beta}$ and 
\begin{align} \left( \asyOp^\alpha_\beta \left( \prod_{l=1}^L (g^{l}(x))^{t_l} \right)\right)(x) = \sum_{l=1}^L t_l (g^l(x))^{t_l-1} \left(\prod_{\substack{m=1 \\ m\neq l}}^L (g^{m}(x))^{t_m}\right) ( \asyOp^\alpha_\beta g^{l})(x) . \end{align}
\end{corollary}
\begin{corollary}
\label{crll:polynomial_comp}
If $g^{1}, \ldots, g^{L} \in \fring{x}{\alpha}{\beta}$  and $p\in \R[y_1,\ldots,y_L]$ is polynomial in $L$ variables, then $p(g^{1}(x), \ldots, g^{L}(x))\in\fring{x}{\alpha}{\beta}$ and
\begin{align} \label{eqn:polynomial_comp} (\asyOp^\alpha_\beta (p(g^{1}, \ldots, g^{L})))(x) = \sum_{l=1}^L \frac{\partial p}{\partial y_l} (g^{1}, \ldots, g^{L}) (\asyOp^\alpha_\beta g^{l})(x). \end{align}
\end{corollary}
Although the last three statements are only basic general properties of commutative derivation rings, they suggest that $\asyOp^\alpha_\beta$ fulfills a simple chain rule. In fact, Corollary~\ref{crll:polynomial_comp} can still be generalized from polynomials to analytic functions (as we will do in Theorem~\ref{thm:chain_analytic}), but, as already mentioned, this intuition turns out to be false in general. 

We will prove Proposition~\ref{prop:derivation} alongside with another statement which will be useful to establish the chain rule:
\begin{proposition}
\label{prop:submultiplicative}
If $\alpha,\beta \in \R_{>0}$ and $R\in \N_0$ are kept fixed, then there exists a constant $C \in \R$ such that 
\begin{align} \label{eqn:submultiplicative} \rho^\alpha_{\beta,R}(f \cdot g) &\leq C \rho^\alpha_{\beta,R}(f) \rho^\alpha_{\beta,R}(g) && \forall f,g \in \fring{x}{\alpha}{\beta}. \end{align}
\end{proposition}
\begin{corollary}
\label{crll:estimate_for_product_pow}
If $\alpha,\beta \in \R_{>0}$, $R\in \N_0$ and $g^{1}, \ldots, g^{L} \in \fring{x}{\alpha}{\beta}$ are kept fixed, then there exists a constant $C \in \R$ such that
\begin{align} \rho^\alpha_{\beta,R} \left( \prod_{l=1}^L (g^{l}(x))^{t_l} \right) &\leq C^{|\textbf{t}|} && \forall \textbf{t} \in \N_0^L \text{ with } |\textbf{t}| \geq 1. \end{align}
\end{corollary}
\begin{proof}
Iterating eq.\ \eqref{eqn:submultiplicative} gives a constant $C\in \R$ such that
\begin{align*} \rho^\alpha_{\beta,R} \left( \prod_{l=1}^L (g^{l}(x))^{t_l} \right) &\leq {C}^{|\textbf{t}|-1} \prod_{l=1}^L \left(\rho^\alpha_{\beta,R} (g^l) \right)^{t_l} &&\forall \textbf{t} \in \N_0^L \text{ with } |\textbf{t}| \geq 1. \end{align*} 
The right hand side is clearly bounded by ${C'}^{|\textbf{t}|}$ for all $|\textbf{t}| \geq 1$ with an appropriate $C'\in \R$ which depends on the $g^l$.
\end{proof}

We will prove Proposition~\ref{prop:derivation} under the assumption that $\beta > 0$. The following lemma shows that, as a consequence of Proposition~\ref{prop:betashift}, we can do so without loss of generality. 
\begin{lemma}
\label{lmm:wlogbetaderivation}
If Proposition~\ref{prop:derivation} holds for $\beta \in \R_{>0}$, then it holds for all $\beta \in \R$.
\end{lemma}
\begin{proof}
For $\beta \in \R$, choose $m\in \N_0$ such that $\beta+m >0$. If $f,g \in \fring{x}{\alpha}{\beta}$, then $f,g \in \fring{x}{\alpha}{\beta+m}$ by Proposition~\ref{prop:betashift}. By the requirement $f\cdot g \in \fring{x}{\alpha}{\beta+m}$ and $( \asyOp^\alpha_{\beta+m} ( f \cdot g ))(x) = f(x) (\asyOp^\alpha_{\beta+m} g)(x) + g(x) (\asyOp^\alpha_{\beta+m} f)(x)$. Using $(\asyOp^\alpha_{\beta+m} f)(x) = x^m (\asyOp^\alpha_{\beta} f)(x)$ from Proposition~\ref{prop:betashift} gives $( \asyOp^\alpha_{\beta+m} ( f \cdot g ))(x) = x^m \left( f(x) (\asyOp^\alpha_{\beta} g)(x) + g(x) (\asyOp^\alpha_{\beta} f)(x) \right)$. Because $f \cdot g \in \fring{x}{\alpha}{\beta+m}$ and $\asyOp^\alpha_{\beta+m} ( f \cdot g ) \in x^m \R[[x]]$, it follows that $f \cdot g \in \fring{x}{\alpha}{\beta}$ and $\asyOp^\alpha_{\beta} ( f \cdot g ) = f(x) (\asyOp^\alpha_{\beta} g)(x) + g(x) (\asyOp^\alpha_{\beta} f)(x)$ by Proposition~\ref{prop:betashift}.
\end{proof}

To prove Propositions~\ref{prop:derivation} and \ref{prop:submultiplicative}, we will use some estimates for the coefficients of the product of two power series. To establish these estimates, we will require that $\beta > 0$.
\begin{lemma}
\label{lmm:estimate_full_fg}
If $\alpha,\beta \in \R_{>0}$ and $R \in \N_0$ are kept fixed, then there exists a constant $C \in \R$ such that for all $f,g \in \fring{x}{\alpha}{\beta}$ and $n,K \in \N_0$ with $K \leq R$ as well as $n \geq K$,
\begin{align} \left| \sum_{m=0}^n f_{n-m} g_m - \sum_{m=0}^{K-1} f_{n-m} g_m - \sum_{m=0}^{K-1} f_m g_{n-m} \right| &\leq C \rho^\alpha_{\beta,R}(f) \rho^\alpha_{\beta,R}(g) \G{n-K}{\alpha}{\beta}. \end{align}
\end{lemma}
\begin{proof} 
Observation~\ref{obs:smallc_bigC_estimate} with $K=0$ states that $\left| f_n \right| \leq \rho^\alpha_{\beta,R}(f) \G{n}{\alpha}{\beta}$ for all $f \in \fring{x}{\alpha}{\beta}$ and $n\in \N_0$.
We can use this to estimate the expression 
\begin{gather*} h_n:=\left| \sum_{m=0}^n f_{n-m} g_m - \sum_{m=0}^{K-1} f_{n-m} g_m - \sum_{m=0}^{K-1} f_m g_{n-m} \right| \end{gather*}
in different ranges for $n$,
\begin{align*} 2K> n \geq K&\Rightarrow & h_n&=&\left| \sum_{m=n-K+1}^{K-1} f_{n-m} g_m\right| & \leq \rho^\alpha_{\beta,R}(f) \rho^\alpha_{\beta,R}(g) \sum_{m=n-K+1}^{K-1} \G{n-m}{\alpha}{\beta} \G{m}{\alpha}{\beta} \\ n \geq 2K&\Rightarrow & h_n&=&\left| \sum_{m=K}^{n-K} f_{n-m} g_m\right| & \leq \rho^\alpha_{\beta,R}(f) \rho^\alpha_{\beta,R}(g) \sum_{m=K}^{n-K} \G{n-m}{\alpha}{\beta} \G{m}{\alpha}{\beta}. \end{align*}
It is trivial to find a constant $C$ such that $\sum_{m=n-K+1}^{K-1} \G{n-m}{\alpha}{\beta} \G{m}{\alpha}{\beta} \leq C \G{n-K}{\alpha}{\beta}$ for all $K\leq R$ and $2K> n \geq K$, because $R$ is fixed and only finitely many inequalities need to be fulfilled.
Corollary~\ref{crll:centersum} guarantees that we can also find a constant $C$ for the second case.
\end{proof}

\begin{lemma}
\label{lmm:estimate_partial_fg}
If $\alpha, \beta \in \R_{>0}$ and $R \in \N_0$ are kept fixed, then there exists a constant $C \in \R$ such that for all $f,g \in \fring{x}{\alpha}{\beta}$ and $n,K \in \N_0$ with $K \leq R$ as well as $n \geq K$,
\begin{align} \left| \sum_{m=0}^{K-1} f_{n-m} g_m - \sum_{k=0}^{K-1} d^{f,g}_k \G{n-k}{\alpha}{\beta} \right| &\leq C \rho^\alpha_{\beta,R}(f) \rho^\alpha_{\beta,R}(g) \G{n-K}{\alpha}{\beta}, \end{align}
where $d^{f,g}_k := [x^k] g(x) (\asyOp^\alpha_\beta f)(x)$.
\end{lemma}
\begin{proof} 
Observation~\ref{obs:smallc_bigC_estimate} with the substitutions $n\rightarrow n-m$ and $K \rightarrow K-m$ implies that
\begin{align*} \left| f_{n-m} - \sum_{k=0}^{K-m-1} c^f_{k} \G{n-m-k}{\alpha}{\beta} \right| \leq \rho^\alpha_{\beta,R}(f) \G{n-K}{\alpha}{\beta}, \end{align*}
for all $f \in \fring{x}{\alpha}{\beta}$ and $n,K,m \in \N_0$ with $m \leq K \leq R$ as well as $n\geq K$ where $c_k^f = [x^k] (\asyOp^\alpha_\beta f)(x)$.  It also follows from Observation~\ref{obs:smallc_bigC_estimate} that $|g_m| \leq \rho^\alpha_{\beta,R}(g) \G{m}{\alpha}{\beta}$ for all $g \in \fring{x}{\alpha}{\beta}$ and $m \in \N_0$.
Because $d^{f,g}_k = \sum_{m=0}^k c^f_{k-m} g_m$,
\begin{gather*} \left| \sum_{m=0}^{K-1} f_{n-m} g_m - \sum_{k=0}^{K-1} d^{f,g}_k \G{n-k}{\alpha}{\beta} \right| = \left| \sum_{m=0}^{K-1} f_{n-m} g_m - \sum_{k=0}^{K-1} \sum_{m=0}^k c^f_{k-m} g_m \G{n-k}{\alpha}{\beta} \right| \\ = \left| \sum_{m=0}^{K-1} \left( f_{n-m} - \sum_{k=m}^{K-1} c^f_{k-m} \G{n-k}{\alpha}{\beta} \right) g_m \right| \leq \sum_{m=0}^{K-1}\left| f_{n-m} - \sum_{k=0}^{K-m-1} c^f_{k} \G{n-m-k}{\alpha}{\beta} \right| \left| g_m \right| \\ \leq \rho^\alpha_{\beta,R}(f) \rho^\alpha_{\beta,R}(g) \G{n-K}{\alpha}{\beta} \sum_{m=0}^{K-1} \G{m}{\alpha}{\beta} \qquad \forall n \geq K.   \end{gather*}
Setting $C = \sum_{m=0}^{R-1} \G{m}{\alpha}{\beta}$ results in the statement.
\end{proof}

\begin{lemma}
\label{lmm:estimate_final_fg}
If $\alpha, \beta \in \R_{>0}$ and $R \in \N_0$ are kept fixed, then there exists a constant $C \in \R$ such that for all $f,g \in \fring{x}{\alpha}{\beta}$ and $n,K \in \N_0$ with $K \leq R$ as well as $n \geq K$,
\begin{align} \label{eqn:estimate_final_fg} \left| \sum_{m=0}^{n} f_{n-m} g_m - \sum_{k=0}^{K-1} c^{f \cdot g}_k \G{n-k}{\alpha}{\beta} \right| &\leq C \rho^\alpha_{\beta,R}(f) \rho^\alpha_{\beta,R}(g) \G{n-K}{\alpha}{\beta}, \end{align}
where $c^{f \cdot g}_k := [x^k] \left( f(x) (\asyOp^\alpha_\beta g)(x) + g(x) (\asyOp^\alpha_\beta f)(x) \right)$.
\end{lemma}
\begin{proof}
Note that $c^{f \cdot g}_k = d^{f,g}_k + d^{g,f}_k$ with $d^{f,g}_k$ from Lemma~\ref{lmm:estimate_partial_fg} and $d^{g,f}_k$ respectively with the roles of $f$ and $g$ switched.
We can use the triangle inequality to deduce that
\begin{gather*} \left| \sum_{m=0}^n f_{n-m} g_m - \sum_{k=0}^{K-1} c^{f \cdot g}_k \G{n-k}{\alpha}{\beta} \right| \leq \left| \sum_{m=0}^n f_{n-m} g_m - \sum_{m=0}^{K-1} f_{n-m} g_m - \sum_{m=0}^{K-1} f_{m} g_{n-m} \right| \\ + \left| \sum_{m=0}^{K-1} f_{n-m} g_m - \sum_{k=0}^{K-1} d^{f,g}_k \G{n-k}{\alpha}{\beta} \right| + \left| \sum_{m=0}^{K-1} f_{m} g_{n-m} - \sum_{k=0}^{K-1} d^{g,f}_k \G{n-k}{\alpha}{\beta} \right|. \end{gather*}
Using Lemmas~\ref{lmm:estimate_full_fg} and \ref{lmm:estimate_partial_fg}  on the respective terms on the right hand side of this inequality results in the statement.
\end{proof}

\begin{proof}[Proof of Proposition~\ref{prop:derivation}]
By Lemma~\ref{lmm:wlogbetaderivation}, it is sufficient to prove Proposition~\ref{prop:derivation} for $\beta > 0$. Therefore, we can apply Lemma~\ref{lmm:estimate_final_fg} for $f,g \in \fring{x}{\alpha}{\beta}$. Eq.\ \eqref{eqn:estimate_final_fg} with $K=R$ directly implies that 
\begin{align*} [x^n] f(x) g(x) = \sum_{m=0}^n f_{n-m} g_m &= \sum_{k=0}^{R-1} c^{f \cdot g}_k \G{n-k}{\alpha}{\beta} + \bigO \left(\G{n-R}{\alpha}{\beta} \right) && \forall R \in \N_0, \end{align*}
with $c^{f \cdot g}_k = [x^k] \left( f(x) (\asyOp^\alpha_\beta g)(x) + g(x) (\asyOp^\alpha_\beta f)(x) \right)$.
By Definition~\ref{def:Fpowerseries}, it follows that $f \cdot g \in \fring{x}{\alpha}{\beta}$ and from Definition~\ref{def:basic_asymp_definition} follows eq.\ \eqref{eqn:leipniz}. 
\end{proof}

\begin{proof}[Proof of Proposition~\ref{prop:submultiplicative}]
If $f,g \in \fring{x}{\alpha}{\beta}$, then $f \cdot g \in \fring{x}{\alpha}{\beta}$ by Proposition~\ref{prop:derivation}.
Because $\beta >0$, we have by Definition~\ref{defn:basic_bigC_estimate}
\begin{align*} \rho^\alpha_{\beta,R}(f \cdot g) &= \max_{0\leq K\leq R} \sup_{n\geq K} \frac{\left| \sum_{m=0}^n f_{n-m} g_m - \sum _{k=0}^{K-1} c_{k}^{f\cdot g} \G{n-k}{\alpha}{\beta} \right|}{\G{n-K}{\alpha}{\beta}} && \forall f,g \in \fring{x}{\alpha}{\beta}, \end{align*}
which is bounded by $C \rho^\alpha_{\beta,R}(f) \rho^\alpha_{\beta,R}(g)$ with some fixed $C \in \R$ as follows directly from Lemma~\ref{lmm:estimate_final_fg}.
\end{proof}
\section{Composition}
\label{sec:comp}

\subsection{Composition by analytic functions}
\begin{theorem}
\label{thm:chain_analytic}
If $\alpha \in \R_{>0}$, $\beta \in \R$, $f \in \R\{y_1,\ldots,y_L\}$ is a function in $L$ variables, which is analytic at the origin, and $g^{1}, \ldots, g^{L} \in \fring{x}{\alpha}{\beta} \cap x\R[[x]]$, then
\begin{itemize}
\item
The composition $f\left(g^{1}(x), \ldots, g^{L}(x)\right)$ is in $\fring{x}{\alpha}{\beta}$.
\item $\asyOp^\alpha_\beta$ fulfills a multivariate chain rule for the composition with analytic functions,
\begin{gather} \label{eqn:chain_analytic} \left(\asyOp^\alpha_\beta f\left(g^{1}, \ldots, g^{L}\right) \right)(x) = \sum_{l=1}^L \frac{\partial f}{\partial y_l} \left(g^{1}, \ldots, g^{L}\right) (\asyOp^\alpha_\beta g^{l})(x). \end{gather}
\end{itemize}
\end{theorem}
In \cite{bender1975asymptotic} Edward Bender established this theorem for the case $L=1$ in a less `generatingfunctionology' based notation. 
If, for example, $g\in \fring{x}{\alpha}{\beta}$ and $f \in \R\{x,y\}$, then his Theorem 1 allows us to calculate the asymptotic expansion of the coefficients of the power series $f(g(x),x)$. 
In fact, Bender analyzed more general power series including series with coefficients which grow even more rapidly than factorially. 

The following proof of Theorem~\ref{thm:chain_analytic} is a straightforward generalization of Bender's Lemma 2 and Theorem 1 in \cite{bender1975asymptotic} to the multivariate case $f \in \R\{y_1,\ldots,y_L\}$.

Again, we will start by verifying that we may assume $\beta > 0$ during the proof of Theorem~\ref{eqn:chain_analytic}.
\begin{lemma}
\label{lmm:wlogbetachainanalytic}
If Theorem~\ref{thm:chain_analytic} holds for $\beta\in \R_{>0}$, then it also holds for all $\beta \in \R$.
\end{lemma}
\begin{proof}
For $\beta\in \R$, choose an $m \in \N_0$ such that $\beta+m > 0$.
If $g^{1}, \ldots, g^{L} \in \fring{x}{\alpha}{\beta} \cap x \R[[x]]$, then by Proposition~\ref{prop:betashift}, $g^{1}, \ldots, g^{L} \in \fring{x}{\alpha}{\beta+m} \cap x \R[[x]]$, 
$(\asyOp^\alpha_{\beta+m} g^{l})(x) = x^m (\asyOp^\alpha_{\beta} g^{l})(x)$ and by the requirement $h(x) := f(g^{1}(x), \ldots, g^{L}(x)) \in \fring{x}{\alpha}{\beta+m}$ as well as 
\begin{gather*} (\asyOp^\alpha_{\beta+m} h)(x) = \sum_{l=1}^L \frac{\partial f}{\partial y_l} (g^{1}, \ldots, g^{L}) (\asyOp^\alpha_{\beta+m} g^{l})(x) = x^m \sum_{l=1}^L \frac{\partial f}{\partial y_l} (g^{1}, \ldots, g^{L}) (\asyOp^\alpha_{\beta} g^{l})(x). \end{gather*}
Due to Proposition~\ref{prop:betashift}, $h\in \fring{x}{\alpha}{\beta}$ and 
$(\asyOp^\alpha_{\beta} h)(x) = \sum_{l=1}^L \frac{\partial f}{\partial y_l} (g^{1}, \ldots, g^{L}) (\asyOp^\alpha_{\beta} g^{l})(x)$.
\end{proof}
As before, we will use our freedom to assume that $\beta>0$ to establish an estimate on the coefficients of products of power series in $x\fring{x}{\alpha}{\beta}$.
\begin{lemma}
\label{lmm:power_prod_estimate}
If $\alpha,\beta \in \R_{>0}$ and $g^{1}, \ldots, g^{L} \in \fring{x}{\alpha}{\beta} \cap x\R[[x]]$ are kept fixed, then there exists a constant $C \in \R$ such that 
\begin{align} \left| [x^n] \prod_{l=1}^L \left( g^l(x) \right)^{t_l} \right| &\leq C^{|\textbf{t}|} \G{n-|\textbf{t}|+1}{\alpha}{\beta} && \forall \textbf{t} \in \N_0^L, n \in \N_0 \text{ with } n \geq |\textbf{t}| \geq 1. \end{align}
\end{lemma}

\begin{proof}

By Proposition~\ref{prop:betashiftlow}, it follows from $g^l \in \fring{x}{\alpha}{\beta} \cap x\R[[x]]$ that $\frac{g^l(x)}{x} \in \fring{x}{\alpha}{\beta+1}$ and therefore by Corollary~\ref{crll:chain_for_product_pow} that $\prod_{l=1}^L \left(\frac{g^{l}(x)}{x}\right)^{t_l} \in \fring{x}{\alpha}{\beta+1}$ for all $\textbf{t} \in \N_0^L$. 
We can apply Observation~\ref{obs:smallc_bigC_estimate} with $R=K=0$ to obtain for all $n \geq |\textbf{t}|$,
\begin{gather*} \left| [x^n] \prod_{l=1}^L \left(g^{l}(x)\right)^{t_l} \right| = \left| [x^{n-|\textbf{t}|}] \prod_{l=1}^L \left(\frac{g^l(x)}{x}\right)^{t_l} \right| \leq \rho^\alpha_{\beta+1,0} \left( \prod_{l=1}^L \left(\frac{g^l(x)}{x}\right)^{t_l} \right)\G{n-|\textbf{t}|}{\alpha}{\beta+1}. \end{gather*}
The statement follows from Corollary~\ref{crll:estimate_for_product_pow} and $\G{n-|\textbf{t}|}{\alpha}{\beta+1} = \G{n-|\textbf{t}|+1}{\alpha}{\beta}$.
\end{proof}

\begin{proof}[Proof of Theorem~\ref{thm:chain_analytic}]
The composition $f(g^{1}(x), \ldots, g^{L}(x))$ can be expressed as the sum
\begin{gather*} \sum_{\substack{\textbf{t} \in \N_0^{L}}} f_{t_1,\ldots,t_L} \prod_{l=1}^L \left(g^{l}(x)\right)^{t_l}, \end{gather*}
which can be split in preparation for the extraction of the coefficients and their asymptotics:
\begin{align*} f(g^{1}(x), \ldots, g^{L}(x))&= \sum_{\substack{\textbf{t} \in \N_0^{L} \\ |\textbf{t}| \leq R}} f_{t_1,\ldots,t_L} \prod_{l=1}^L \left(g^{l}(x)\right)^{t_l} + \sum_{\substack{\textbf{t} \in \N_0^{L} \\ |\textbf{t}| > R}} f_{t_1,\ldots,t_L} \prod_{l=1}^L \left(g^{l}(x)\right)^{t_l} && \forall R\in \N_0 . \end{align*}
The first sum is just the composition by a polynomial. 
Corollary~\ref{crll:polynomial_comp} asserts that this sum is in $\fring{x}{\alpha}{\beta}$. It has the asymptotic expansion given in eq.\ \eqref{eqn:polynomial_comp} which agrees, as a series in $x$, with the right hand side of eq.\ \eqref{eqn:chain_analytic} up to order $R-1$, because the partial derivative reduces the order of a polynomial by one and $g^l_0=0$.

It is left to prove that the coefficients of the power series given by the remaining sum over $|\textbf{t}|>R$ are in $\bigO(\G{n-R}{\alpha}{\beta})$. 
Because of Lemma~\ref{lmm:wlogbetachainanalytic}, we may assume that $\beta > 0$ without loss of generality and apply Lemma~\ref{lmm:power_prod_estimate}. Together with the fact that there is a constant $C$, such that $|f_{t_1,\ldots,t_L}|\leq {C}^{|\textbf{t}|}$ for all $\textbf{t} \in \N_0^L$, due to the analyticity of $f$, Lemma~\ref{lmm:power_prod_estimate} ensures that there is a constant $C'\in \R$ such that
\begin{gather*}   \left|[x^n]\sum_{\substack{\textbf{t} \in \N_0^{L} \\ |\textbf{t}| > R}} f_{t_1,\ldots,t_L} \prod_{l=1}^L \left(g^{l}(x)\right)^{t_l}\right| \leq \sum_{\substack{\textbf{t} \in \N_0^{L} \\ n \geq |\textbf{t}| > R}} \left| f_{t_1,\ldots,t_L} \right| \left|[x^n] \prod_{l=1}^L \left(g^{l}(x)\right)^{t_l} \right|     \\ \leq \sum_{t=R+1}^n {C'}^{t} \G{n-t+1}{\alpha}{\beta} \sum_{\substack{\textbf{t} \in \N_0^{L} \\ |\textbf{t}| = t }} 1      ,      \end{gather*}
for all $n \geq R+1$. Because the last sum $| \{ t_1,\ldots, t_L \in \N_0 | t_1 + \cdots + t_L = t \}| = { t + L -1 \choose L-1 }$ is a polynomial in $t$,
Corollary~\ref{crll:exponential_inO} asserts that this is in $\bigO\left(\G{n-R}{\alpha}{\beta}\right)$. 
\end{proof}
\subsection{Proof of the main theorem: Composition of power series in $\fring{x}{\alpha}{\beta}$}
Despite the fact that Bender's theorem applies to a broader range of compositions $f \circ g$, where $f$ does not need to be analytic and $g$ does not need to be an element of $\fring{x}{\alpha}{\beta}$, it cannot be used in the case $f,g \in \fring{x}{\alpha}{\beta}$, where $f \notin \ker \asyOp^\alpha_\beta$. The problem is that we cannot truncate the sum $\sum_{k=0}^\infty f_k g(x)^k$ without losing significant information. In this section, we will confront this problem and prove the general chain rule for the asymptotic derivative:
\begin{theorem}
\label{thm:chainrule}
If $\alpha \in \R_{>0}$, $\beta \in \R$ and $f,g \in \fring{x}{\alpha}{\beta}$ with $g_0=0$ and $g_1=1$, then
\begin{itemize}
\item
The composition $f \circ g$ and the inverse $g^{-1}$ belong to $\fring{x}{\alpha}{\beta}$.
\item 
$\asyOp^\alpha_\beta$ fulfills a chain rule and there is a formula for the $\asyOp^\alpha_\beta$-derivative of the compositional inverse:
\begin{align} \label{eqn:asymp_chainrule_normal} (\asyOp^\alpha_\beta (f \circ g)) (x) = f'(g(x)) (\asyOp^\alpha_\beta g)(x) + \left(\frac{x}{g(x)} \right)^{\beta} e^{\frac{\frac{1}{x} - \frac{1}{g(x)}}{\alpha}} (\asyOp^\alpha_\beta f ) (g(x)), \\ \label{eqn:inverse_asymp} (\asyOp^\alpha_\beta g^{-1})(x) = -{(g^{-1})}'(x) \left(\frac{x}{g^{-1}(x)}\right)^\beta e^{\frac{\frac{1}{x} - \frac{1}{g^{-1}(x)}}{\alpha }} (\asyOp^\alpha_\beta g ) (g^{-1}(x)). \end{align}
\end{itemize}
\end{theorem}

\begin{corollary}
\label{crll:chainimplicit}
If $f \in \R[[x]]$, $g \in \fring{x}{\alpha}{\beta}$ with $g_0=0$, $g_1=1$ and $f\circ g \in \fring{x}{\alpha}{\beta}$, then $f \in \fring{x}{\alpha}{\beta}$.
\end{corollary}
\begin{proof}
Theorem~\ref{thm:chainrule} guarantees that $g^{-1} \in \fring{x}{\alpha}{\beta}$ and therefore also $f = (f \circ g) \circ g^{-1} \in \fring{x}{\alpha}{\beta}$. 
\end{proof}

As before, we will assume that $\beta >0$ while proving Theorem~\ref{thm:chainrule}. The following lemma establishes that we can do so.

\begin{lemma}
\label{lmm:wloggeneralchainrolebeta}
If Theorem~\ref{thm:chainrule} holds for $\beta \in \R_{>0}$, then it holds for all $\beta\in \R$.
\end{lemma}
\begin{proof}
For $\beta \in \R$, choose $m\in \N_0$ such that $\beta+m >0$. 
If $f,g \in \fring{x}{\alpha}{\beta}$ with $g_0=0,g_1=1$, then $f,g \in \fring{x}{\alpha}{\beta+m}$ by Proposition~\ref{prop:betashift}. Because of $(\asyOp^\alpha_{\beta+m} f ) (x) = x^m(\asyOp^\alpha_{\beta} f ) (x)$ and by the requirement
\begin{gather*} (\asyOp^\alpha_{\beta+m} (f \circ g)) (x) = f'(g(x)) (\asyOp^\alpha_{\beta+m} g)(x) + \left(\frac{x}{g(x)} \right)^{\beta+m} e^{\frac{\frac{1}{x} - \frac{1}{g(x)}}{\alpha}} (\asyOp^\alpha_{\beta+m} f ) (g(x)) \\ = x^m \left( f'(g(x)) (\asyOp^\alpha_{\beta} g)(x) + \left(\frac{x}{g(x)} \right)^{\beta} e^{\frac{\frac{1}{x} - \frac{1}{g(x)}}{\alpha}} (\asyOp^\alpha_{\beta} f ) (g(x)) \right). \end{gather*}
Applying Proposition~\ref{prop:betashift} again results in $f \circ g \in \fring{x}{\alpha}{\beta}$ and eq. \eqref{eqn:asymp_chainrule_normal}. Eq. \eqref{eqn:inverse_asymp} and $g^{-1} \in \fring{x}{\alpha}{\beta}$ follow analogously. 
\end{proof}

Obviously, $x \in \fring{x}{\alpha}{\beta}$. We will use this basic fact to prove Theorem~\ref{thm:chainrule} by ensuring that from $f,g \in \fring{x}{\alpha}{\beta}$ follows $f \circ g^{-1} \in \fring{x}{\alpha}{\beta}$ and by constructing the asymptotic expansion of the coefficients of $(f \circ g^{-1})(x)$. To prove that $f \circ g^{-1} \in \fring{x}{\alpha}{\beta}$, the Lagrange inversion formula, which involves the usual derivative of a power series, will be used. To handle this derivative, the following proposition will become necessary:

\begin{proposition}
\label{prop:ordinary_derivative}
If $f\in \fring{x}{\alpha}{\beta}$, then $f'(x) \in \fring{x}{\alpha}{{\beta+2}}$ and 
\begin{align} \label{eqn:derivative_commute} (\asyOp^\alpha_{\beta+2} f')(x) =\left(\alpha^{-1} - x \beta + x^2 \frac{\partial}{\partial x}\right)(\asyOp^\alpha_\beta f)(x). \end{align}
\end{proposition}
\begin{proof}
Recall that $f'(x) = \sum_{n=0}^\infty n f_n x^{n-1} = \sum_{n=0}^\infty (n+1)f_{n+1} x^n$. If $f\in \fring{x}{\alpha}{\beta}$, then by Definition~\ref{def:Fpowerseries}, 
\begin{align*} (n+1)f_{n+1} &= \sum _{k=0}^{R-1} c^f_k (n+1) \G{n+1-k}{\alpha}{\beta} + (n+1)\bigO\left(\G{n+1-R}{\alpha}{\beta}\right) && \forall R\in \N_0. \end{align*}
Observe that because $x \Gamma(x) = \Gamma(x+1)$ and $\G{n}{\alpha}{\beta}=\alpha^{n+\beta}\Gamma(n+\beta)$,
\begin{gather*} (n+1) \G{n+1-k}{\alpha}{\beta} \\ = \alpha^{n+1-k+\beta} \left( (n+1-k+\beta) \Gamma(n+1-k+\beta) + (k-\beta) \Gamma(n+1-k+\beta) \right) \\ = \alpha^{-1} \G{n-k}{\alpha}{\beta+2} + (k-\beta) \G{n-k-1}{\alpha}{\beta+2}. \end{gather*}
Therefore, for all $R\in \N_0$
\begin{gather*} (n+1)f_{n+1} = \sum _{k=0}^{R-1} c^f_k \left( \alpha^{-1} \G{n-k}{\alpha}{\beta+2} + (k-\beta) \G{n-k-1}{\alpha}{\beta+2}\right) + \bigO\left(\G{n-R}{\alpha}{\beta+2}\right), \end{gather*}
and it follows from Definition~\ref{def:Fpowerseries} that $f'\in \fring{x}{\alpha}{{\beta+2}}$.
Moreover, by Definition~\ref{def:basic_asymp_definition}, 
\begin{gather*} (\asyOp^\alpha_{\beta+2} f')(x) = \sum_{k=0}^\infty c^{f'}_k x^k = \sum_{k=0}^\infty c^f_k \left( \alpha^{-1} x^k + (k-\beta) x^{k+1}\right) \\ = \left(\alpha^{-1} - x \beta + x^2 \frac{\partial}{\partial x}\right) (\asyOp^\alpha_\beta f)(x). \qedhere \end{gather*}
\end{proof}

While using the Lagrange inversion formula to establish $f \circ g^{-1} \in \fring{x}{\alpha}{\beta}$, it will be convenient to work in the rings $\fring{x}{\alpha}{\beta+1}$ and $\fring{x}{\alpha}{\beta+2}$, which contain $\fring{x}{\alpha}{\beta}$ as a subring. Therefore, we will start with some observations on intermediate quantities in $\fring{x}{\alpha}{\beta+1}$ and $\fring{x}{\alpha}{\beta+2}$. The following three lemmas are basic applications of the chain rule for the composition with analytic functions and the product rule, but we will prove them in detail to get acquainted to the new notions from the last sections.

\begin{lemma}
    \label{lmm:chainasympterm0}
    If $g \in \fring{x}{\alpha}{\beta}$ with $g_0=0,g_1=1$ and $\gamma \in \R$, then $\left(\frac{g(x)}{x}\right)^{\gamma} \in \fring{x}{\alpha}{\beta+1}$ and
    \begin{align} \left( \asyOp^\alpha_{\beta+1} \left(\frac{g(x)}{x}\right)^{\gamma} \right) = \gamma \left(\frac{g(x)}{x}\right)^{\gamma-1} \left(\asyOp^\alpha_{\beta}g\right)(x). \end{align}
\end{lemma}
\begin{proof}
    Observe that $F(x) := (1-x)^\gamma \in \R\{x\}$ and $F'(x) = -\gamma (1-x)^{\gamma-1}$. Proposition~\ref{prop:betashiftlow} implies that $\frac{g(x)}{x} \in \fring{x}{\alpha}{\beta+1}$, because $g\in \fring{x}{\alpha}{\beta} \cap x\R[[x]]$. As $g_1 = 1$, we additionally have $1- \frac{g(x)}{x} \in \fring{x}{\alpha}{\beta+1} \cap x\R[[x]]$. Using Theorem~\ref{thm:chain_analytic} results in
    \begin{gather*} \left(\frac{g(x)}{x}\right)^{\gamma} = F\left( 1- \frac{g(x)}{x} \right) \in \fring{x}{\alpha}{\beta+1}, \intertext{and by the chain rule for the composition with analytic functions from eq.\ \eqref{eqn:chain_analytic},} \left( \asyOp^\alpha_{\beta+1} \left(\frac{g(x)}{x}\right)^{\gamma} \right) = F'\left( 1- \frac{g(x)}{x} \right) \left(\asyOp^\alpha_{\beta+1} \left( 1- \frac{g(x)}{x}\right) \right)(x) \\ = -\gamma \left(\frac{g(x)}{x}\right)^{\gamma-1}\left( \asyOp^\alpha_{\beta+1}\left( -\frac{g(x)}{x}\right) \right)(x) = \gamma \left(\frac{g(x)}{x}\right)^{\gamma-1} \left(\asyOp^\alpha_{\beta}g\right)(x), \end{gather*}
    where we used the linearity of $\asyOp^\alpha_{\beta+1}$ and $\left( \asyOp^\alpha_{\beta+1}\frac{g(x)}{x} \right)(x) = \left( \asyOp^\alpha_{\beta} g \right)(x)$ due to Proposition~\ref{prop:betashiftlow}.
\end{proof}
\begin{lemma}
\label{lmm:chainasympinterterm1}
    If $g \in \fring{x}{\alpha}{\beta}$ with $g_0=0,g_1=1$, then 
    \begin{gather} \label{eqn:chainasympinterterm1infring} A(x) := \frac{1}{g(x)}-\frac{1}{x} \in \fring{x}{\alpha}{\beta+2}\text{, } e^{\frac{A(x)}{\alpha}} \in \fring{x}{\alpha}{\beta+2} \text{ and } \\ \label{eqn:chainasympinterterm1} \left( \asyOp^\alpha_{\beta+2} e^{\frac{A(x)}{\alpha}} \right)(x) = -\alpha^{-1} \left( \frac{x}{g(x)}\right)^2 e^{\frac{A(x)}{\alpha}} \left( \asyOp^\alpha_{\beta} g \right)(x). \end{gather}
\end{lemma}
\begin{proof}
From Lemma~\ref{lmm:chainasympterm0} with $\gamma=-1$, it follows that $\frac{x}{g(x)} \in \fring{x}{\alpha}{\beta+1}$ and 
\begin{gather*} \left( \asyOp^\alpha_{\beta+1} \frac{x}{g(x)} \right) = - \left(\frac{x}{g(x)}\right)^{2} \left(\asyOp^\alpha_{\beta}g\right)(x). \end{gather*}
Because $g_1=1$, $\frac{x}{g(x)}-1 \in \fring{x}{\alpha}{\beta+1} \cap x\R[[x]]$. Moreover, by Proposition~\ref{prop:betashiftlow}, $A(x) = \frac{\frac{x}{g(x)}-1}{x} \in \fring{x}{\alpha}{\beta+2}$ and 
\begin{gather} \label{eqn:derivativeb2A} \left( \asyOp^\alpha_{\beta+2} A \right) (x) = \left( \asyOp^\alpha_{\beta+1} \left( \frac{x}{g(x)}-1\right) \right)(x) = - \left(\frac{x}{g(x)}\right)^{2} \left(\asyOp^\alpha_{\beta}g\right)(x). \end{gather}
Observe that $\frac{A(x)-A(0)}{\alpha} \in \fring{x}{\alpha}{\beta+2}\cap x \R[[x]]$. Because $e^x \in \R\{x\}$, we can apply Theorem~\ref{thm:chain_analytic} to conclude that $e^{\frac{A(x)-A(0)}{\alpha}} \in \fring{x}{\alpha}{\beta+2}$ and by linearity that also $e^{\frac{A(x)}{\alpha}} \in \fring{x}{\alpha}{\beta+2}$. Finally, we can use the chain rule for the composition with analytic functions to write the left hand side of eq.\ \eqref{eqn:chainasympinterterm1} as
\begin{gather*} e^{\frac{A(0)}{\alpha}} \left( \asyOp^\alpha_{\beta+2} e^{ \frac{A(x)-A(0)}{\alpha} } \right) (x) = e^{\frac{A(0)}{\alpha}} e^{ \frac{A(x)-A(0)}{\alpha} } \left( \asyOp^\alpha_{\beta+2} \frac{A(x)-A(0)}{\alpha} \right) (x) \\ = e^{\frac{A(x)}{\alpha} } \left( \asyOp^\alpha_{\beta+2} \frac{A(x)}{\alpha} \right) (x) = \alpha^{-1} e^{\frac{A(x)}{\alpha} } \left( \asyOp^\alpha_{\beta+2} A \right) (x). \end{gather*}
The statement in eq.\ \eqref{eqn:chainasympinterterm1} follows after substitution of $\left( \asyOp^\alpha_{\beta+2} A \right) (x)$ from eq.\ \eqref{eqn:derivativeb2A}.
\end{proof}

\begin{lemma}
\label{lmm:chainasympinterterm2}
If $f, g \in \fring{x}{\alpha}{\beta}$ with $g_0=0,g_1=1$ and $\gamma \in \R$, then 
\begin{gather} \label{eqn:chainasympinterterm2infring} B_\gamma(x) := f(x) g'(x) \left(\frac{g(x)}{x}\right)^{\gamma} \in \fring{x}{\alpha}{\beta+2} \qquad \text{ and } \\ \label{eqn:chainasympinterterm2} \begin{gathered} \left( \asyOp^\alpha_{\beta+2} B_\gamma \right)(x) = \\ \left(\frac{g(x)}{x}\right)^{\gamma} \left( x^2 g'(x)\left( \asyOp^\alpha_{\beta}f \right)(x) + f(x) \left( \gamma x g'(x) \frac{x}{g(x)} + \alpha^{-1} - \beta x + x^2 \frac{\partial}{\partial x} \right)(\asyOp^\alpha_{\beta} g)(x) \right). \end{gathered} \end{gather}
\end{lemma}
\begin{proof}
Recall that due to Proposition~\ref{prop:betashift}, $f\in \fring{x}{\alpha}{\beta} \subset \fring{x}{\alpha}{\beta+2}$ and $(\asyOp^\alpha_{\beta+2} f)(x) = x^2 (\asyOp^\alpha_{\beta} f)(x)$. Proposition~\ref{prop:ordinary_derivative} guarantees that $g' \in \fring{x}{\alpha}{\beta+2}$ and 
\begin{align*} (\asyOp^\alpha_{\beta+2} g')(x) &f= \left(\alpha^{-1} - x \beta + x^2 \frac{\partial}{\partial x}\right)(\asyOp^\alpha_\beta g)(x). \end{align*}
Because of Lemma~\ref{lmm:chainasympterm0} and Proposition~\ref{prop:betashift}, we have $\left(\frac{g(x)}{x}\right)^{\gamma}\in \fring{x}{\alpha}{\beta+2}$ and 
\begin{align*} \left(\asyOp^\alpha_{\beta+2} \left(\frac{g(x)}{x}\right)^{\gamma}\right)(x) = x \gamma \left(\frac{g(x)}{x}\right)^{\gamma-1} \left(\asyOp^\alpha_{\beta}g\right)(x). \end{align*}
Putting all this together we can use Corollary~\ref{crll:chain_for_product} with $g^1(x) = f(x)$, $g^2(x)= g'(x)$ and $g^3(x) = \left(\frac{g(x)}{x}\right)^{\gamma}$ to obtain eqs.\ \eqref{eqn:chainasympinterterm2infring} and \eqref{eqn:chainasympinterterm2}.
\end{proof}

\begin{lemma}
\label{lmm:ABrho}
If $\alpha,\beta \in \R_{>0}$, $R \in \N_0$ and $A,B_\gamma$ as defined in eqs.\ \eqref{eqn:chainasympinterterm1infring} and \eqref{eqn:chainasympinterterm2infring} are kept fixed, then there exists a constant $C \in \R$ such that
\begin{align} \rho^\alpha_{\beta+2,R}\left( B_\gamma(x) A(x)^m \right) &\leq C^{m+1} && \forall m\in \N_0. \end{align}
\end{lemma}
\begin{proof}
Apply Corollary~\ref{crll:chain_for_product_pow} with $g^1(x) = B_\gamma(x), g^2(x) =A(x), t_1 = 1$ and $t_2 =m$ to verify that $B_\gamma(x) A(x)^m \in \fring{x}{\alpha}{\beta+2}$ for all $m\in \N_0$. Apply Corollary~\ref{crll:estimate_for_product_pow} with the same parameters.
\end{proof}
\begin{corollary}
\label{crll:ABestimate}
If $\alpha,\beta \in \R_{>0}$, $R \in \N_0$ and $A,B_\gamma$ as defined in eqs.\ \eqref{eqn:chainasympinterterm1infring} and \eqref{eqn:chainasympinterterm2infring} are kept fixed, then there exists a constant $C \in \R$ such that 
\begin{align*} \left| [x^{n}] B_\gamma(x) A(x)^m - \sum_{k=0}^{R-1} c_{k,m} \G{n-k}{\alpha}{\beta+2} \right| &\leq C^{m+1} \G{n-R}{\alpha}{\beta+2} && \forall n \geq R \text{ and } m\in \N_0 \end{align*}
where $c_{k,m} = [x^k] \left( \asyOp^\alpha_{\beta+2} B_\gamma(x) A(x)^m \right)(x)$.
\end{corollary}
\begin{proof}
Additionally to Lemma~\ref{lmm:ABrho}, apply Observation~\ref{obs:smallc_bigC_estimate} with $K=R$.
\end{proof}

The key to the extraction of the large $n$ asymptotics of $[x^n] (f \circ g^{-1})(x)$ is a variant of the Chu-Vandermonde identity. We will prove this identity using elementary power series techniques.

\begin{lemma}
\label{lmm:identity}
For all $a\in \R$ and $m,k\in \N_0$
\begin{align} \label{eqn:chuidentity} { a \choose m } &= \sum_{l=0}^m {k+l-1 \choose l} { a -k-l \choose m-l }. \end{align}
\end{lemma}
\begin{proof}
Recall that ${ a \choose n }= [x^n] (1+x)^a$ for all $a\in \R$ and $n\in \N_0$. 
By standard generating function arguments it follows that $[x^n] \frac{1}{(1-x)^k} = { k + n - 1 \choose n }$ for all $n,k\in \N_0$.
Observe that for all $a\in \R$ and $k\in \N_0$, we have the following identities in $\R[[x]]$:
\begin{gather*} (1+x)^{a} = (1+x)^k (1+x)^{a-k} = \frac{1}{\left( 1-\frac{x}{1+x}\right)^k}(1+x)^{a-k} \\ = \sum_{l=0}^\infty { k + l - 1 \choose l }\left(\frac{x}{1+x}\right)^l (1+x)^{a-k} = \sum_{l=0}^\infty { k + l - 1 \choose l } x^l (1+x)^{a-k-l}. \end{gather*}
Extracting coefficients from the first and the last expression results in the Chu-Vander\-monde-type identity in eq.\ \eqref{eqn:chuidentity}.
\end{proof}

\begin{corollary}
\label{crll:identity_rdy}
For all $\alpha,\beta \in \R_{>0}$ and $n,R,k \in \N_0$ with $n\geq R\geq k$, we have the identity in $\R[x]$
\begin{gather} \sum_{m=0}^{n-R} {n+\beta+1 \choose m} \G{n-m-k}{\alpha}{\beta+2} x^m = \sum \limits_{l = 0}^{n-R} { l+k-1 \choose l } \G{n-l-k}{\alpha}{\beta+2} x^{l} \sum_{m=0}^{n-R-l} \frac{\left(\frac{x}{\alpha}\right)^{m}}{m!}. \end{gather}
\end{corollary}
\begin{proof}
Observe that ${ a \choose n } = \frac{1}{n!} \frac{\Gamma(a+1)}{\Gamma(a-n+1)}$ for all $a\in \R$ and $n\in \N_0$ as long as $n< a+1$.
By writing the second binomial coefficient on the right hand side of eq.\ \eqref{eqn:chuidentity} in this form and setting $a= n+\beta+1$, we get for all 
$n,m,k \in \N_0$ with 
$m+k < n+\beta+2$
\begin{align*} { n+\beta+1 \choose m } \Gamma( n - m - k+\beta+2) = \sum \limits_{l = 0}^m { k+l-1 \choose l } \frac{\Gamma(n-k-l+\beta+2)}{(m-l)!}. \end{align*}
Multiplying by $x^m \alpha^{n-m-k+\beta+2}$, summing over $m$ and using $\G{n}{\alpha}{\beta} = \alpha^{n+\beta} \Gamma(n+\beta)$ gives,
\begin{align*} \sum_{m=0}^{n-R} { n+\beta+1 \choose m } \G{ n - m - k }{\alpha}{\beta+2} x^m = \sum_{m=0}^{n-R} \sum \limits_{l = 0}^m { k+l-1 \choose l } \frac{\alpha^{l-m} \G{n-k-l}{\alpha}{\beta+2}}{(m-l)!} x^m . \end{align*}
Note that $k \leq R$ and $m \leq n-R$ imply $m+k \leq n < n+ \beta+2$. The statement follows after changing the order of summation and a shift of the summation variable $m\rightarrow m+l$ both on the right hand side of this equation.
\end{proof}

We are now equipped with the necessary tools to tackle the asymptotic analysis of the coefficients of $(f \circ g^{-1})(x)$. The first step is to express $(f \circ g^{-1})(x)$ in terms of the intermediate power series $A(x)$ and $B_\gamma(x)$. We will do so using a variant of the Lagrange inversion theorem.
\begin{lemma}
\label{lmm:lagrange}
If $p,q \in \R[[x]]$ with $q_0=0$ and $q_1 = 1$, then
\begin{align} [x^n] p\left( q^{-1}(x)\right) &= [x^n] p(x) q'(x) \left(\frac{x}{q(x)}\right)^{n+1} && \forall n\in \N_0. \end{align}
\end{lemma}
\begin{proof}
Note that the identity holds for $n=0$, because $q_0=0$ and $q_1 = 1$. It follows from the Lagrange inversion theorem \cite[A.6]{flajolet2009analytic} for $n\geq 1$ that 
\begin{gather*} [x^n] p\left(q^{-1}(x)\right) = \frac{1}{n}[x^{n-1}] p'(x) \left(\frac{x}{q(x)}\right)^{n} \\ = \frac{1}{n}[x^{n-1}] \frac{\partial}{\partial x} \left( p(x) \left(\frac{x}{q(x)}\right)^{n} \right) - \frac{1}{n}[x^{n-1}] p(x) \left( \frac{\partial}{\partial x} \left(\frac{x}{q(x)}\right)^{n} \right).    \end{gather*}
Using $\frac{1}{n}[x^{n-1}] \frac{\partial}{\partial x} = [x^n]$ and evaluating the derivative in the second term result in the statement.
\end{proof}

\begin{corollary}
\label{crll:chainrule_lagrange}
If $\alpha,\beta \in \R_{>0}$, $f, g \in \fring{x}{\alpha}{\beta}$ and $A, B_\gamma$ as defined in eqs.\ \eqref{eqn:chainasympinterterm1infring} and \eqref{eqn:chainasympinterterm2infring}, then
\begin{align} \label{eqn:diffgroupsumrep1} [x^n] f(g^{-1}(x)) &= \sum_{m=0}^n {n+\beta+1 \choose m} [x^{n-m}] B_\beta(x) A(x)^m && \forall n\in \N_0. \end{align}
\end{corollary}

\begin{proof}
By Lemma~\ref{lmm:lagrange},
\begin{gather*} [x^n] f(g^{-1}(x)) = [x^{n}] f(x)g'(x)\left(\frac{x}{g(x)}\right)^{n+1} = [x^{n}] f(x)g'(x)\left(\frac{g(x)}{x}\right)^{\beta} \left(\frac{x}{g(x)}\right)^{n+\beta+1}. \end{gather*}
Using the definitions of $A$ and $B_\gamma$ gives $[x^n] f(g^{-1}(x)) = [x^{n}] B_\beta(x) \left(1+x A(x)\right)^{n+\beta+1}$. Expanding with the generalized binomial theorem results in eq.\ \eqref{eqn:diffgroupsumrep1}.
\end{proof}

\begin{corollary}
\label{crll:diffgroupsumsubstituted_1}
If $\alpha,\beta \in \R_{>0}$, $f, g \in \fring{x}{\alpha}{\beta}$ and $A, B_\gamma$ as defined in eqs.\ \eqref{eqn:chainasympinterterm1infring} and \eqref{eqn:chainasympinterterm2infring}, then
\begin{align} \label{eqn:diffgroupsumsubstituted_1} [x^n] f(g^{-1}(x)) &= \sum_{m=0}^{n-R} {n+\beta+1 \choose m} [x^{n-m}] B_\beta(x) A(x)^m + \bigO\left(\G{n-R}{\alpha}{\beta+2}\right)&& \forall R\in \N_0. \end{align}
\end{corollary}
\begin{proof}
Eq.\ \eqref{eqn:diffgroupsumsubstituted_1} follows from eq.\ \eqref{eqn:diffgroupsumrep1} and 
\begin{gather*} \left| \sum_{m=n-R+1}^{n} {n+\beta+1 \choose m} [x^{n-m}] B_\beta(x) A(x)^{m} \right| = \left|\sum_{m=0}^{R-1} {n+\beta+1 \choose n-m} [x^{m}] B_\beta(x) A(x)^{n-m} \right| \\ \leq \sum_{m=0}^{R-1} {n+\beta+1 \choose n-m} C^{n-m+1} \G{m}{\alpha}{\beta+2} \in \bigO\left(\G{n-R}{\alpha}{\beta+2}\right) \qquad \forall R\in \N_0, \end{gather*}
where the second step, together with the existence of an appropriate $C\in \R$, follows from Corollary~\ref{crll:ABestimate} with $R=0$ and the inclusion holds, because ${n+\beta+1 \choose n-m} = \frac{\Gamma(n+\beta+2)}{\Gamma(n-m+1) \Gamma(\beta+m+2)} \sim \frac{n^{\beta+m+1}}{\Gamma(\beta+m+2)}$ by elementary properties of the $\Gamma$ function.
\end{proof}

\begin{lemma}
\label{lmm:diffgroupsumsubstituted_3}
If $\alpha,\beta \in \R_{>0}$, $f, g \in \fring{x}{\alpha}{\beta}$ and $A, B_\gamma$ as defined in eqs.\ \eqref{eqn:chainasympinterterm1infring} and \eqref{eqn:chainasympinterterm2infring}, then 
\begin{multline}
\label{eqn:diffgroupsumsubstituted_3}
[x^n] f(g^{-1}(x)) = \sum_{k=0}^{R-1}\sum_{l = 0}^{n-R} \sum_{m=0}^{n-R-l} c_{k,l,m} { l+k-1 \choose l } \G{n-l-k}{\alpha}{\beta+2} 
+ \bigO\left(\G{n-R}{\alpha}{\beta+2}\right)
\\
\forall R\in \N_0,
\end{multline}
where $c_{k,l,m} := [x^k] \left( \asyOp^\alpha_{\beta+2} B_\beta(x) A(x)^{l} \frac{\left(\frac{A(x)}{\alpha}\right)^{m}}{m!} \right)(x)$. 
\end{lemma}
Note that the terms of the triple sum in eq.\ \eqref{eqn:diffgroupsumsubstituted_3} where $k=0$ are not all trivial, because ${-1 \choose 0}=1$ by the definition of the binomial coefficients with negative arguments.

\begin{proof}

For all $n,m \in \N_0$ with $n-m \geq R$ set
\begin{align*} \mathcal{R}_{n,m}:=[x^{n-m}] B_\beta(x) A(x)^m - \sum_{k=0}^{R-1} c_{k,m} \G{n-m-k}{\alpha}{\beta+2}, \end{align*}
where $c_{k,m} = [x^k] \left( \asyOp^\alpha_{\beta+2} B_\beta(x) A(x)^m \right)(x)$.
Substituting $\mathcal{R}_{n,m}$ into eq.\ \eqref{eqn:diffgroupsumsubstituted_1} gives
\begin{gather} \begin{gathered} \label{eqn:diffgroupsumsubstituted_2} [x^n] f(g^{-1}(x)) = \sum_{m=0}^{n-R} {n+\beta+1 \choose m} \sum_{k=0}^{R-1} c_{k,m} \G{n-m-k}{\alpha}{\beta+2} \\ + \sum_{m=0}^{n-R} {n+\beta+1 \choose m} \mathcal{R}_{n,m} + \bigO\left(\G{n-R}{\alpha}{\beta+2}\right) \qquad \forall R\in \N_0, \end{gathered} \end{gather}
By Corollary~\ref{crll:ABestimate} with $n\rightarrow n-m$, we can find a constant $C \in \R$ such that 
\begin{align*} \left| \mathcal{R}_{n,m} \right| &\leq C^{m+1} \G{n-m-R}{\alpha}{\beta+2} && \forall n-m \geq R. \intertext{Therefore,} \mathcal{R}_n :=\left| \sum_{m=0}^{n-R} {n+\beta+1 \choose m} \mathcal{R}_{n,m} \right| &\leq \sum_{m=0}^{n-R} {n+\beta+1 \choose m} C^{m+1} \G{n-m-R}{\alpha}{\beta+2} &&\forall n \geq R. \end{align*}
Applying Corollary~\ref{crll:identity_rdy} with $x\rightarrow C$ and $k=R$ results in
\begin{align*} \mathcal{R}_n &\leq C \sum \limits_{l = 0}^{n-R} { l+R-1 \choose l } \G{n-l-R}{\alpha}{\beta+2} C^{l} \sum_{m=0}^{n-R-l} \frac{\left(\frac{C}{\alpha}\right)^{m}}{m!} && \forall n\geq R \\ &\leq C \sum \limits_{l = R}^{n} { l-1 \choose l-R } \G{n-l}{\alpha}{\beta+2} C^{l-R} \sum_{m=0}^{n-R} \frac{\left(\frac{C}{\alpha}\right)^{m}}{m!} && \forall n\geq R. \end{align*}
From $\sum_{m=0}^{n-R} \frac{\left(\frac{C}{\alpha}\right)^{m}}{m!} \leq e^{\frac{C}{\alpha}}$ and Corollary~\ref{crll:exponential_inO}, it follows that $\mathcal{R}_{n} \in \bigO\left(\G{n-R}{\alpha}{\beta+2}\right)$, because ${ l-1 \choose l-R }$ is a polynomial in $l$. Therefore, for all $R\in \N_0$
\begin{gather*} [x^n] f(g^{-1}(x)) = \sum_{m=0}^{n-R} {n+\beta+1 \choose m} \sum_{k=0}^{R-1} c_{k,m} \G{n-m-k}{\alpha}{\beta+2} + \bigO\left(\G{n-R}{\alpha}{\beta+2}\right) \\ = \sum_{k=0}^{R-1} [x^k] \left( \asyOp^\alpha_{\beta+2} B_\beta(x) \sum_{m=0}^{n-R} {n+\beta+1 \choose m} A(x)^m \G{n-m-k}{\alpha}{\beta+2} \right) + \bigO\left(\G{n-R}{\alpha}{\beta+2}\right) , \end{gather*}
where $\asyOp^\alpha_{\beta+2}$ acts on everything on its right. Applying Corollary~\ref{crll:identity_rdy} with $x \rightarrow A(x)$ to the inner sum and reordering result in the statement.
\end{proof}
\begin{lemma}
\label{lmm:diffgroupsumsubstituted_5}
If $\alpha,\beta \in \R_{>0}$, $R \in \N_0$ and $A,B_\gamma$ as defined in eqs.\ \eqref{eqn:chainasympinterterm1infring} and \eqref{eqn:chainasympinterterm2infring}, then 
\begin{align} \label{eqn:diffgroupsumsubstituted_5} [x^n] f(g^{-1}(x)) &= \sum_{k=0}^{R-1}\sum_{l = 0}^{R-1-k} c_{k,l}' { l+k-1 \choose l } \G{n-l-k}{\alpha}{\beta+2} + \bigO\left(\G{n-R}{\alpha}{\beta+2}\right) &&\forall R \in \N_0, \end{align}
where $c_{k,l}' := [x^k] \left( \asyOp^\alpha_{\beta+2} B_\beta(x) A(x)^{l} e^{ \frac{A(x)}{\alpha} } \right)(x)$.
\end{lemma}

\begin{proof}
Set $c_{k,l,m}$ as in Lemma~\ref{lmm:diffgroupsumsubstituted_3}.
By Lemma~\ref{lmm:ABrho} there exists a constant $C \in \R$ such that 
$\rho^\alpha_{\beta+2,R}\left( B_\beta(x) A(x)^{l+m} \right) \leq C^{l+m+1}$ for all $l,m\in \N_0$. It follows from the second part of Observation~\ref{obs:smallc_bigC_estimate} that
\begin{align} \label{eqn:casympestimate} |c_{k,l,m}| &= \frac{\alpha^{-m}}{m!} \left| [x^k]\left( \asyOp^\alpha_{\beta+2} B_\beta(x) A(x)^{l+m} \right)(x) \right| \leq \frac{\alpha^{-m}}{m!} C^{l+m+1} && \forall k,l,m\in \N_0 \text{ with } k \leq R. \end{align}
Therefore, for all $k\leq R$ and $n \geq 2R - k$,
\begin{gather*} \left| \sum_{l = R-k}^{n-R} \sum_{m=0}^{n-R-l} c_{k,l,m} { l+k-1 \choose l } \G{n-l-k}{\alpha}{\beta+2} \right| \\ \leq \sum_{l = R-k}^{n-R} \sum_{m=0}^{n-R-l} \frac{\alpha^{-m}C^{l+m+1}}{m!} { l+k-1 \choose l } \G{n-l-k}{\alpha}{\beta+2} \end{gather*}
which is in $\bigO\left(\G{n-R}{\alpha}{\beta+2}\right)$, because $\sum_{m=0}^{n-R-l}\frac{\alpha^{-m}C^{m}}{m!} \leq e^{\frac{C}{\alpha}}$ and by Corollary~\ref{crll:exponential_inO}. Applying this to truncate the summation over $l$ in eq.\ \eqref{eqn:diffgroupsumsubstituted_3} from Lemma~\ref{lmm:diffgroupsumsubstituted_3} gives for all $R \in \N_0$
\begin{align} \label{eqn:diffgroupsumsubstituted_4} [x^n] f(g^{-1}(x)) &= \sum_{k=0}^{R-1}\sum_{l = 0}^{R-k-1} \sum_{m=0}^{n-R-l} c_{k,l,m} { l+k-1 \choose l } \G{n-l-k}{\alpha}{\beta+2} + \bigO\left(\G{n-R}{\alpha}{\beta+2}\right). \end{align}
Note that ${ n + m \choose n} \geq 1 \Rightarrow (n+m)! \geq n! m!$ and therefore
\begin{gather*} \sum_{m=n}^\infty \frac{C^m}{m!} = \sum_{m=0}^\infty \frac{C^{n+m}}{(n+m)!} \leq \frac{C^n}{n!} \sum_{m=0}^\infty \frac{C^{m}}{m!} = e^C \frac{C^n}{n!}. \end{gather*}
It follows from this and eq.\ \eqref{eqn:casympestimate} that for all $n\geq R-l+1$ and $k+l\leq R$
\begin{gather*} \left| \sum_{m=n-R-l+1}^{\infty} c_{k,l,m} \G{n-l-k}{\alpha}{\beta+2}\right| \leq C^{l+1} \sum_{m=n-R-l+1}^{\infty} \frac{\left(\frac{C}{\alpha}\right)^m}{m!} \G{n-l-k}{\alpha}{\beta+2} \\ \leq e^{\frac{C}{\alpha}} C^{l+1} \left(\frac{C}{\alpha}\right)^{n-l-R+1} \frac{\G{n-l-k}{\alpha}{\beta+2}}{(n-R-l+1)!}, \end{gather*}
which is in $\bigO\left(\G{n-R}{\alpha}{\beta+2}\right)$ as long as $k$ and $l$ are bounded, because  $\frac{\Gamma(n-l-k+\beta+2)}{\Gamma(n-R-l+2)} \sim n^{R-k+\beta}$. Applying this to complete the summation over $m$ in eq.\ \eqref{eqn:diffgroupsumsubstituted_4} and noting that $c_{k,l}' = \sum_{m=0}^\infty c_{k,l,m}$ results in eq.\ \eqref{eqn:diffgroupsumsubstituted_5}.
\end{proof}
\begin{corollary}
\label{crll:chainrule_asymp_expansion_1}
If $\alpha,\beta \in \R_{>0}$, $R \in \N_0$ and $A,B_\gamma$ as defined in eqs.\ \eqref{eqn:chainasympinterterm1infring} and \eqref{eqn:chainasympinterterm2infring}, then $f \circ g^{-1} \in \fring{x}{\alpha}{\beta+2}$ and 
\begin{align} \label{eqn:chainrule_asymp_expansion_1} [x^k] \left( \asyOp^\alpha_{\beta+2} f \circ g^{-1} \right)(x) &= [x^{k}] \left( \asyOp^\alpha_{\beta+2} B_{\beta-k+1}(x) e^{\frac{A(x)}{\alpha} }\right)(x) && \forall k\in \N_0. \end{align}
\end{corollary}
\begin{proof}
After the change of summation variables $k \rightarrow k+l$, eq. \eqref{eqn:diffgroupsumsubstituted_5} becomes
\begin{gather*} [x^n] f(g^{-1}(x)) = \sum_{k=0}^{R-1}\sum_{l = 0}^{k} c_{k-l,l}' { k-1 \choose l } \G{n-k}{\alpha}{\beta+2} + \bigO\left(\G{n-R}{\alpha}{\beta+2}\right) \qquad\forall R \in \N_0. \end{gather*}
By Definition~\ref{def:Fpowerseries}, this equation states that $f \circ g^{-1} \in \fring{x}{\alpha}{\beta+2}$ and that 
the coefficients of the asymptotic expansion are 
\begin{gather*} c^{f \circ g^{-1}}_k = \sum_{l = 0}^{k} c_{k-l,l}' { k-1 \choose l } = \sum_{l = 0}^{k} [x^{k-l}] \left( \asyOp^\alpha_{\beta+2} B_\beta(x) A(x)^{l} { k-1 \choose l } e^{\frac{A(x)}{\alpha}} \right)(x) \\ =[x^{k}] \sum_{l = 0}^{\infty} x^l \left( \asyOp^\alpha_{\beta+2} B_\beta(x) A(x)^{l} { k-1 \choose l } e^{\frac{A(x)}{\alpha}} \right)(x) \\ = [x^{k}] \left( \asyOp^\alpha_{\beta+2} B_\beta(x) \sum_{l=0}^\infty (x A(x))^{l} { k-1 \choose l } e^{\frac{A(x)}{\alpha}} \right)(x), \end{gather*}
where $x^l \left( \asyOp^\alpha_{\beta+2} f(x)\right)(x) = \left( \asyOp^\alpha_{\beta+2} x^l f(x)\right)(x)$ for all $f\in \fring{x}{\alpha}{\beta+2}$ was used, which follows from the product rule (Proposition~\ref{prop:derivation}). Because of $\sum_{l=0}^\infty { k-1 \choose l } (x A(x))^{l} = (1+xA(x))^{k-1} = \left( \frac{x}{g(x)} \right)^{k-1}$ and the definition of $B_\gamma$ in Lemma~\ref{lmm:chainasympinterterm2}, the statement follows.
\end{proof}

\begin{proof}[Proof of Theorem~\ref{thm:chainrule}]
Because of Lemma~\ref{lmm:wloggeneralchainrolebeta}, we may assume that $\beta \in \R_{>0}$ and start with the expression from Corollary~\ref{crll:chainrule_asymp_expansion_1} for $[x^k] \left( \asyOp^\alpha_{\beta+2} f \circ g^{-1} \right)(x)$. We will use Lemmas~\ref{lmm:chainasympinterterm1} and \ref{lmm:chainasympinterterm2} to expand this expression. By Corollary~\ref{crll:chainrule_asymp_expansion_1} and the product rule (Proposition~\ref{prop:derivation}), we have for all $k\in \N_0$
\begin{align} \label{eqn:thereallyfirstandsecondterm} [x^k] \left( \asyOp^\alpha_{\beta+2} f \circ g^{-1} \right)(x) &= [x^{k}] \left( e^{\frac{A(x)}{\alpha} }\left(\asyOp^\alpha_{\beta+2} B_{\beta-k+1} \right) (x) + B_{\beta-k+1}(x) \left( \asyOp^\alpha_{\beta+2} e^{\frac{A(x)}{\alpha} }\right)(x) \right). \end{align}
Applying Lemma~\ref{lmm:chainasympinterterm2} on the first term of this expression gives after a straightforward but lengthy calculation,
\begin{gather} \begin{gathered} \label{eqn:thereallyfirstterm} [x^{k}] e^{\frac{A(x)}{\alpha} }\left(\asyOp^\alpha_{\beta+2} B_{\beta-k+1} \right) (x) = [x^{k}] e^{\frac{A(x)}{\alpha} }\left(\frac{g(x)}{x}\right)^{\beta-k+1} \Bigg( x^2 g'(x) \left( \asyOp^\alpha_{\beta}f \right)(x) \\ + f(x) \left( x(\beta-k+1) g'(x) \frac{x}{g(x)} + \alpha^{-1} - \beta x + x^2 \frac{\partial}{\partial x} \right)(\asyOp^\alpha_{\beta} g)(x) \Bigg) \\ = [x^{k}] e^{\frac{A(x)}{\alpha} }\left(\frac{g(x)}{x}\right)^{\beta-k+1} \Bigg( x^2 g'(x) \left( \asyOp^\alpha_{\beta}f \right)(x) \\ +\left(- x^2 f'(x) + \alpha^{-1}f(x) g'(x) \left( \frac{x}{g(x)} \right)^2 \right) (\asyOp^\alpha_{\beta} g)(x) \Bigg), \end{gathered} \end{gather}
where the identity $[x^k] x p'(x) q(x) = k [x^k] p(x) q(x) - [x^k] x p(x) q'(x)$ for all $p,q \in \R[[x]]$ was used to eliminate the summand which contains the $\frac{\partial}{\partial x} (\asyOp^\alpha_{\beta} g)(x)$ factor.
By Lemma~\ref{lmm:chainasympinterterm1}, the second term on the right hand side of eq.\ \eqref{eqn:thereallyfirstandsecondterm} is 
\begin{gather} \begin{gathered} \label{eqn:thereallysecondterm} [x^k] B_{\beta-k+1}(x) \left( \asyOp^\alpha_{\beta+2} e^{\frac{A(x)}{\alpha} }\right)(x) = -[x^k] \alpha^{-1} B_{\beta-k+1}(x) \left( \frac{x}{g(x)}\right)^2 e^{\frac{A(x)}{\alpha}} \left( \asyOp^\alpha_{\beta} g \right)(x) \\ = -[x^k]\alpha^{-1} f(x) g'(x) \left(\frac{g(x)}{x}\right)^{\beta-k+1} \left( \frac{x}{g(x)}\right)^2 e^{\frac{A(x)}{\alpha}} \left( \asyOp^\alpha_{\beta} g \right)(x), \end{gathered} \end{gather}
where the definition of $B_{\beta-k+1}(x)$ from Lemma~\ref{lmm:chainasympinterterm2} was substituted. Summing both expressions for the terms in eq.\ \eqref{eqn:thereallyfirstandsecondterm} from eqs.\ \eqref{eqn:thereallyfirstterm} and \eqref{eqn:thereallysecondterm} and substituting the definition of $A(x)$ from Lemma~\ref{lmm:chainasympinterterm1} results in 
\begin{align*} [x^k] \left( \asyOp^\alpha_{\beta+2} f \circ g^{-1} \right)(x)&= [x^k]x^2 e^{\frac{\frac{1}{g(x)}-\frac{1}{x}}{\alpha} } \left(\frac{g(x)}{x}\right)^{\beta-k+1} \left( g'(x)(\asyOp^\alpha_{\beta} f)(x) - f'(x) (\asyOp^\alpha_{\beta} g)(x) \right), \end{align*}
for all $k\in \N_0$.
By Proposition~\ref{prop:betashift}, the $x^2$ prefactor indicates that $f \circ g^{-1}$ is actually in the subspace $\fring{x}{\alpha}{\beta} \subset \fring{x}{\alpha}{\beta+2}$ and
\begin{align*} [x^k] \left( \asyOp^\alpha_{\beta} f \circ g^{-1} \right)(x)&= [x^k]e^{\frac{\frac{1}{g(x)}-\frac{1}{x}}{\alpha} } \left(\frac{g(x)}{x}\right)^{\beta-k-1} \left( g'(x)(\asyOp^\alpha_{\beta} f)(x) - f'(x) (\asyOp^\alpha_{\beta} g)(x) \right). \end{align*}
If we set $p(x) := e^{\frac{\frac{1}{g(x)}-\frac{1}{x}}{\alpha} }\left(\frac{g(x)}{x}\right)^{\beta} \left( (\asyOp^\alpha_{\beta} f)(x) - \frac{f'(x)}{g'(x)} (\asyOp^\alpha_{\beta} g)(x) \right)$ and $q(x) := g(x)$, we obtain
\begin{align*} [x^k] \left( \asyOp^\alpha_{\beta} f \circ g^{-1} \right)(x) &= [x^{k}] p(x) q'(x) \left(\frac{x}{q(x)}\right)^{k+1} = [x^k]p(q^{-1}(x)) && \forall k \in \N_0, \end{align*}
by Lemma~\ref{lmm:lagrange}.
After replacing $p$ and $q$ by their expressions, we obtain 
\begin{align} \label{eqn:conjugate_asymp_equation} (\asyOp^\alpha_{\beta} f \circ g^{-1} )(x) &= e^{\frac{\frac{1}{x}-\frac{1}{g^{-1}(x)}}{\alpha} } \left(\frac{x}{g^{-1}(x)}\right)^{\beta} \left( (\asyOp^\alpha_{\beta} f)(g^{-1}(x)) -\frac{f'(g^{-1}(x))}{g'(g^{-1}(x))} (\asyOp^\alpha_{\beta} g)(g^{-1}(x)) \right). \end{align}
The special case $f(x)=x$ with an application of the identity $g'(g^{-1}(x)) = \frac{1}{{(g^{-1})}'(x)}$ results in eq.\ \eqref{eqn:inverse_asymp}. Solving eq.\ \eqref{eqn:inverse_asymp} for $(\asyOp^\alpha_{\beta} g)(g^{-1}(x))$ and substituting the result into eq.\ \eqref{eqn:conjugate_asymp_equation} gives eq.\ \eqref{eqn:asymp_chainrule_normal} with the substitution $g \rightarrow g^{-1}$.
\end{proof}

\begin{remark}
Bender and Richmond \cite{Bender1984} established that $[x^n](1+g(x))^{\gamma n + \delta} = n \gamma e^{\frac{\gamma g_1}{\alpha}} g_n + \bigO(g_n)$ if $g_n \sim \alpha n g_{n-1}$ and $g_0 = 0$. Using Lagrange inversion, the first coefficient in the expansion of the compositional inverse in eq.\ \eqref{eqn:inverse_asymp} can be obtained from this. In this respect, Theorem~\ref{thm:chainrule} is a generalization of Bender and Richmond's result.

In the same article Bender and Richmond proved a theorem similar to Theorem~\ref{thm:chainrule} for the class of power series $f$ whose coefficients grow more rapidly than factorially such that $n f_{n-1} \in \smallO(f_n)$. Theorem~\ref{thm:chainrule} establishes a link to the excluded case $n f_{n-1} = \bigO(f_n)$.
\end{remark}
\begin{remark}
    The restriction $g_1=1$ ensures that our power series actually have compositional inverses and that we do not leave the ring $\fring{x}{\alpha}{\beta}$. We might also allow a non-zero positive value for $g_1$. To do this it is sufficient to allow composition with the family of power series $h_\gamma(x) = \gamma x \in \R[[x]]$ where $\gamma \in \R_{> 0}$. Right composition with $h_\gamma$, $f \mapsto f \circ h_\gamma$ is a trivial isomorphisms of vector spaces $\fring{x}{\alpha}{\beta} \rightarrow \fring{x}{\gamma \alpha}{\beta}$, which follows immediately from Definition~\ref{def:Fpowerseries}. Every power series $g(x) \in \fring{x}{\alpha}{\beta}$ with $g_0=0$ and $g_1 > 0$ can be decomposed into $g = \widetilde g \circ h_{g_1}$ such that $\widetilde g_1 = 1$ and $\widetilde g \in \fring{x}{\alpha / g_1}{\beta}$. The asymptotics of the coefficients of the composition $f \circ g= f \circ \widetilde g \circ h_{g_1}$ can be calculated using Theorem~\ref{thm:chainrule} if $f,\widetilde g \in \fring{x}{\alpha'}{\beta}$ for some $\alpha' \geq \frac{\alpha}{g_1}$.
\end{remark}

\begin{remark}
The chain rule in eq.\ \eqref{eqn:asymp_chainrule_normal}
exposes a peculiar algebraic structure. It would be useful to have a combinatorial interpretation of the $e^{\frac{\frac{1}{x} - \frac{1}{g(x)}}{\alpha}}$ term. 
\end{remark}
\section{Some remarks on differential equations}
\label{sec:dgls}

Differential equations arising from physical systems form an active field of research in the scope of resurgence \cite{garoufalidis2012asymptotics,aniceto2011resurgence}. Unfortunately, the exact calculation of an overall factor of the asymptotic expansion of a solution of an ODE, called \textit{Stokes constant}, turns out to be difficult for many problems. This fact severely limits the utility of the method for enumeration problems, as the dominant factor of the asymptotic expansion is of most interest and the detailed structure of the asymptotic expansion is secondary.

In this section it will be sketched, for the sake of completeness, how the presented combinatorial framework fits into the realm of differential equations. The given elementary properties each have their counterpart in resurgence's alien calculus \cite[II.6]{mitschi2016divergent}.

Corollary~\ref{crll:polynomial_comp} serves as a good starting point to analyze differential equations with power series solutions in $\fring{x}{\alpha}{\beta}$. Given a polynomial $F \in \R[x, y_0, \ldots,y_L]$, the $\asyOp^\alpha_\beta$-derivation can be applied to the ordinary differential equation
\begin{align*} 0 = F(x, f(x), f'(x), f''(x), \ldots, f^{(L)}(x)). \end{align*}
Applying the $\asyOp$-derivation naively to both sides of this equation and using the chain rule for the composition with polynomials results in a linear equation for the asymptotic expansions of the derivatives $f^{(l)}$. Proposition~\ref{prop:ordinary_derivative} tells us, how the asymptotic expansions of the $f^{(l)}$ relate to each other. 
We will follow this line of thought in detail in
\begin{proposition}
\label{prop:dgl_in_F}
If $F \in \R[x, y_0, \ldots,y_L]$ and $f \in \fring{x}{\alpha}{\beta}$ is a solution of the differential equation
\begin{gather} \label{eqn:initialdgl} 0 = F\left(x, f(x), f'(x), f''(x), \ldots, f^{(L)}(x)\right), \intertext{then $(\asyOp^\alpha_\beta f)(x)$ is a solution of the linear differential equation} \label{eqn:conddgllinasymp} 0 = \sum_{l=0}^L x^{2L-2l} \frac{\partial F}{\partial y_l} \left(x, f^{(0)},\ldots,f^{(L)}\right) \left( \prod_{j=0}^{l-1} \left(\alpha^{-1} - x (\beta + 2j) + x^2 \frac{\partial}{\partial x}\right) \right) (\asyOp^\alpha_\beta f)(x). \end{gather}
\end{proposition}
\begin{proof}
From Proposition~\ref{prop:ordinary_derivative}, Proposition~\ref{prop:betashift} and $f \in \fring{x}{\alpha}{\beta}$, it follows that $f^{(l)} \in \fring{x}{\alpha}{\beta+2l} \subset \fring{x}{\alpha}{\beta+2L}$ for all $L \geq l$.
By Corollary~\ref{crll:polynomial_comp}, we can apply $\asyOp^\alpha_{\beta+2L}$ to both sides of eq.\ \eqref{eqn:initialdgl} and use Proposition~\ref{prop:betashift},
\begin{align} \begin{split} \label{eqn:dgl_in_F_intermediate} 0 &= \sum_{l=0}^L \frac{\partial F}{\partial y_l} \left(x, f^{(0)},\ldots,f^{(L)}\right) \left(\asyOp^\alpha_{\beta+2L} f^{(l)}\right)(x) \\ &= \sum_{l=0}^L \frac{\partial F}{\partial y_l} \left(x, f^{(0)},\ldots,f^{(L)}\right) x^{2(L-l)}\left(\asyOp^\alpha_{\beta+2l} f^{(l)}\right)(x). \end{split} \end{align}
Iterating Proposition~\ref{prop:ordinary_derivative} gives
\begin{align*} \left(\asyOp^\alpha_{\beta+2l} f^{(l)}\right)(x) &= \left(\alpha^{-1} - x (\beta + 2(l-1)) + x^2 \frac{\partial}{\partial x}\right) \left(\asyOp^\alpha_{\beta+2(l-1)} f^{(l-1)}\right)(x) \\ &= \left( \prod_{j=0}^{l-1} \left(\alpha^{-1} - x (\beta + 2j) + x^2 \frac{\partial}{\partial x}\right) \right) (\asyOp^\alpha_\beta f)(x). \end{align*}
Substituting this into eq.\ \eqref{eqn:dgl_in_F_intermediate} results in eq.\ \eqref{eqn:conddgllinasymp}.
\end{proof}

\begin{remark}
Even if it is known that the coefficients of the power series solution of a differential equation have a well-behaved asymptotic expansion, Proposition~\ref{prop:dgl_in_F} provides this asymptotic expansion only up to the initial values for the linear differential equation \eqref{eqn:conddgllinasymp}. Note that the form of the asymptotic expansion can still depend non-trivially on the initial values of the solution $f$ of a non-linear differential equation.
\end{remark}
\begin{remark}
    The linear differential equation \eqref{eqn:conddgllinasymp} only has a non-trivial solution in $\R[[x]]$ if $\alpha^{-1}$ is the root of a certain polynomial. More specifically, making a power series ansatz for $(\asyOp^\alpha_\beta f)(x)$ in eq.\ \eqref{eqn:conddgllinasymp} gives
\begin{gather*} 0 = [x^m]\sum_{l=0}^L x^{2L-2l} \alpha^{-l} \frac{\partial F}{\partial y_l} \left(x, f^{(0)},\ldots,f^{(L)}\right), \end{gather*}
where $m$ is the smallest integer such that the equation is not trivially fulfilled.
If this root is not real or if two roots have the same modulus, the present formalism has to be generalized to complex and multiple $\alpha$ to express the asymptotic expansion of a general solution. This generalization is straightforward. We merely need to generalize Definition~\ref{def:Fpowerseries} of suitable sequences to: 
\end{remark}
\begin{definition}
For given $\beta \in \R$ and $\alpha_1, \ldots, \alpha_L \in \C$ with $|\alpha_1| = |\alpha_2| = \ldots = |\alpha_L| =: \alpha > 0$ let $\cring{x}{\alpha_1, \ldots, \alpha_L}{\beta} \subset \C[[x]]$ be the subspace of complex power series, such that $f \in \cring{x}{\alpha_1, \ldots, \alpha_L}{\beta} $  if and only if there exist sequences of complex numbers $(c_{k,l}^f)_{k\in \N_0, l \in [1,L]}$, which fulfill
\begin{align} f_n &= \sum _{k=0}^{R-1} \sum_{l=1}^L c_{k,l}^f \G{n-k}{{\alpha_l}}{\beta} + \bigO\left(\G{n-k}{{\alpha}}{\beta}\right) && \forall R \in \N_0. \end{align}
\end{definition}

\section{Applications}
\label{sec:applications}

\subsection{Connected chord diagrams}
\begin{figure}%
\begin{subfigure}[b]{0.5\textwidth}%
\centering
\begin{tikzpicture}[scale=0.6] \draw[draw=none,fill=red!20!white] (4.700000,-0.150000) rectangle (6.300000,0.750000); \draw[draw=none,fill=red!20!white] (0.600000,-0.250000) rectangle (4.400000,1.750000); \draw (0,0) arc (0:-180:-3.500000); \draw[color=red!20!white,line width=3pt] (1,0) arc (0:-180:-1.000000); \draw (1,0) arc (0:-180:-1.000000); \draw[color=red!20!white,line width=3pt] (2,0) arc (0:-180:-1.000000); \draw (2,0) arc (0:-180:-1.000000); \draw[color=red!20!white,line width=3pt] (5,0) arc (0:-180:-0.500000); \draw (5,0) arc (0:-180:-0.500000); \node at (0, -.75){$1$}; \draw[fill] (0, 0) circle (1pt); \node at (1, -.75){$2$}; \draw[fill] (1, 0) circle (1pt); \node at (2, -.75){$3$}; \draw[fill] (2, 0) circle (1pt); \node at (3, -.75){$4$}; \draw[fill] (3, 0) circle (1pt); \node at (4, -.75){$5$}; \draw[fill] (4, 0) circle (1pt); \node at (5, -.75){$6$}; \draw[fill] (5, 0) circle (1pt); \node at (6, -.75){$7$}; \draw[fill] (6, 0) circle (1pt); \node at (7, -.75){$8$}; \draw[fill] (7, 0) circle (1pt); \draw (-0.500000,0)--(7.500000,0); \end{tikzpicture}
\subcaption{disconnected chord diagram}
\end{subfigure}%
\begin{subfigure}[b]{0.5\textwidth}%
\centering
\begin{tikzpicture}[scale=0.6] \draw[color=white,line width=3pt] (0,0) arc (0:-180:-2.500000); \draw (0,0) arc (0:-180:-2.500000); \draw[color=white,line width=3pt] (2,0) arc (0:-180:-1.000000); \draw (2,0) arc (0:-180:-1.000000); \draw[color=white,line width=3pt] (1,0) arc (0:-180:-3.000000); \draw (1,0) arc (0:-180:-3.000000); \draw[color=white,line width=3pt] (3,0) arc (0:-180:-1.500000); \draw (3,0) arc (0:-180:-1.500000); \node at (0, -.75){$1$}; \draw[fill] (0, 0) circle (1pt); \node at (1, -.75){$2$}; \draw[fill] (1, 0) circle (1pt); \node at (2, -.75){$3$}; \draw[fill] (2, 0) circle (1pt); \node at (3, -.75){$4$}; \draw[fill] (3, 0) circle (1pt); \node at (4, -.75){$5$}; \draw[fill] (4, 0) circle (1pt); \node at (5, -.75){$6$}; \draw[fill] (5, 0) circle (1pt); \node at (6, -.75){$7$}; \draw[fill] (6, 0) circle (1pt); \node at (7, -.75){$8$}; \draw[fill] (7, 0) circle (1pt); \draw (-0.500000,0)--(7.500000,0); \end{tikzpicture}
\subcaption{connected chord diagram}
\end{subfigure}%
\caption{Illustrations of connected and disconnected chord diagrams. The red rectangles indicate the connected components of the disconnected diagram.}%
\label{fig:disc_cntd_chord}
\end{figure}
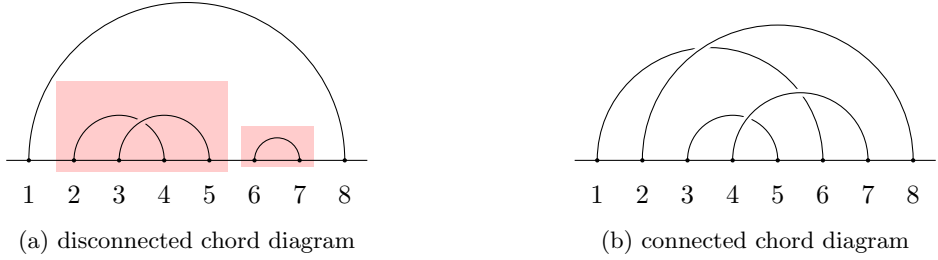

A chord diagram with $n$ chords is a circle with $2n$ points, which are labeled by integers $1,\ldots,2n$ and connected in disjoint pairs by $n$ chords. There are $(2n-1)!!$ such diagrams. 

A chord diagram is \textit{connected} if no set of chords can be separated from the remaining chords by a line which does not cross any chords. The difference between connected and disconnected chord diagrams is illustrated in Figure~\ref{fig:disc_cntd_chord}. Let $I(x) = \sum_{n=0}(2n-1)!!x^n$, the ordinary generating function of \textit{all} chord diagrams, and $C(x) = \sum_{n=0} C_n x^n$, where $C_n$ is the number of \textit{connected} chord diagrams with $n$ chords. Following \cite{Flajolet2000}, the power series $I(x)$ and $C(x)$ are related by,
\begin{align} \label{eqn:cntd_chords} I(x) &= 1+C\left(xI(x)^2\right). \end{align}
This functional equation can be solved for the coefficients of $C(x)$ by basic iterative methods. The first coefficients are 
\begin{align} C(x) = x + x^2 + 4x^3 + 27x^4 + 248 x^5 +\cdots    \end{align}
This sequence is entry \texttt{A000699} in Neil Sloane's integer sequence on-line encyclopedia \cite{oeis}.

Because $(2n-1)!!=\frac{2^{n+\frac12}}{\sqrt{2 \pi}}\Gamma(n+\frac12)=\frac{1}{\sqrt{2\pi}}\G{n}{2}{\frac12}$, the power series $I$ is in $\fring{x}{2}{{\frac12}}$ and $\left(\asyOp^2_{\frac12} I\right)(x) = \frac{1}{\sqrt{2 \pi}}$ as a direct consequence of Definitions~\ref{def:Fpowerseries} and \ref{def:basic_asymp_definition}. From eq.\ \eqref{eqn:cntd_chords}, it also follows that $C(xI(x)^2) \in \fring{x}{2}{{\frac12}}$. Because $x I(x)^2 \in \fring{x}{2}{{\frac12}}$ by the product rule (Proposition~\ref{prop:derivation}), we know from Corollary~\ref{crll:chainimplicit} with $f(x) = C(x)$ and $g(x) = x I(x)^2$ that $C \in \fring{x}{2}{{\frac12}}$.

Applications of the general chain rule from Theorem~\ref{thm:chainrule} and the product rule on the functional eq. \eqref{eqn:cntd_chords} result in
\begin{gather} \begin{gathered}  \left(\asyOp^2_{\frac12} I\right)(x) = \left( \asyOp^2_{\frac12} \left(1+C\left(xI(x)^2\right)\right) \right)(x) = \left( \asyOp^2_{\frac12} C\left(xI(x)^2\right) \right)(x) \\= 2 x I(x) C'\left(xI(x)^2\right) (\asyOp^2_{\frac12} I)(x) + \left(\frac{x}{xI(x)^2} \right)^\frac12 e^{\frac{xI(x)^2-x}{2 x^2 I(x)^2}} \left(\asyOp^2_{\frac12} C\right)\left(xI(x)^2\right). \end{gathered} \end{gather}
which can be solved for $\left(\asyOp^2_{\frac12} C\right)\left(xI(x)^2\right)$,
\begin{align*} \left(\asyOp^2_{\frac12} C\right)\left(xI(x)^2\right) &= \frac{I(x)-2 x I(x)^2 C'\left(xI(x)^2\right)}{\sqrt{2\pi}} e^{\frac{1-I(x)^2}{2 xI(x)^2}}, \end{align*}
where $\left(\asyOp^2_{\frac12} I\right)(x) = \frac{1}{\sqrt{2 \pi}}$ was used. This can be composed with the unique $y \in \R[[x]]$ which solves $y(x) I(y(x))^2 = x$,
\begin{align*} \left(\asyOp^2_{\frac12} C\right)(x) &= \frac{I(y(x))-2 x C'(x)}{\sqrt{2\pi}} e^{\frac{1-I(y(x))^2}{2 x}}. \end{align*}
From eq.\ \eqref{eqn:cntd_chords}, it follows that $I(y(x)) = 1+C(x)$, therefore 
\begin{align}      \label{eqn:asyCpre} \left(\asyOp^2_{\frac12} C\right)(x)&= \frac{1+C(x)-2x C'(x)}{\sqrt{2\pi}}e^{-\frac{1}{2 x}( 2C(x) + C(x)^2 )}. \end{align}
It can be verified, using the closed form of its coefficients, that the power series $I(x)$ fulfills the differential equation $2x^2 I'(x) + xI(x)+1 = I(x)$. From this and eq.\ \eqref{eqn:cntd_chords}, the non-linear differential equation $C'(x) = \frac{C(x)(1+C(x))-x}{2x C(x)}$ \cite{Flajolet2000} for $C(x)$ can be deduced. Using this on the expression for $(\asyOp^2_{\frac12} C)(x)$ from eq.\ \eqref{eqn:asyCpre} results in the simplification,
\begin{align} \label{eqn:asyC} \left(\asyOp^2_{\frac12} C\right)(x)&= \frac{1}{\sqrt{2\pi}}\frac{x}{C(x)}e^{-\frac{1}{2 x}(2C(x) + C(x)^2)}. \end{align}
This is the generating function of the full asymptotic expansion of $C_n$. 
The first coefficients are,
\begin{align} \left(\asyOp^2_{\frac12} C\right)(x) = \frac{e^{-1}}{\sqrt{2 \pi}} \left( 1 - \frac{5}{2} x - \frac{43}{8} x^2 - \frac{579}{16} x^3 - \frac{44477}{128} x^4 - \frac{5326191}{1280} x^5 + \cdots \right).   \end{align}
By Definitions~\ref{def:Fpowerseries} and \ref{def:basic_asymp_definition} as well as $\frac{1}{\sqrt{2\pi}}\G{n}{2}{\frac12} = (2n-1)!!$, we get the two equivalent expressions for the asymptotic expansion of the coefficients $C_n$:
\begin{align*} C_n &=\sum_{k = 0}^{R-1} \G{n -k}{2}{\frac12} [x^k] \left(\asyOp^2_{\frac12} C\right)(x) + \bigO\left(\G{n -R}{2}{\frac12}\right) && \forall R \in \N_0 \\ C_n&= \sqrt{2 \pi} \sum_{k = 0}^{R-1} (2(n-k)-1)!! [x^k] \left(\asyOp^2_{\frac12} C\right)(x) + \bigO\left(\left(2(n-R)-1\right)!!\right)&& \forall R \in \N_0. \end{align*}
The first terms of this large $n$ expansion are
\begin{align*} C_n&= e^{-1} \left( (2n-1)!! - \frac{5}{2} (2n-3)!! - \frac{43}{8} (2n-5)!!    - \frac{579}{16} (2n-7)!!   + \cdots \right).       \end{align*}
The first term, $e^{-1}$, of this expansion has been computed by Kleitman \cite{kleitman1970proportions}, Stein and Everett \cite{stein1978class} and Bender and Richmond \cite{Bender1984} each using different methods. With the presented method an arbitrary number of coefficients can be computed.
Some additional coefficients are given in Table \ref{tab:coefficientsCM}.

\begin{table}
\centering
\tiny{
\def\arraystretch{1.5}
\begin{tabular}{|c||c||c|c|c|c|c|c|c|c|}
\hline
sequence&$0$&$1$&$2$&$3$&$4$&$5$&$6$&$7$&$8$\\
\hline\hline
$e\sqrt{2\pi}(\asyOp^2_{\frac12} C)$&$1$&$- \frac{5}{2}$&$- \frac{43}{8}$&$- \frac{579}{16}$&$- \frac{44477}{128}$&$- \frac{5326191}{1280}$&$- \frac{180306541}{3072}$&$- \frac{203331297947}{215040}$&$- \frac{58726239094693}{3440640}$\\
\hline
$e\sqrt{2\pi}(\asyOp^2_{\frac12} M)$&$1$&$-4$&$-6$&$- \frac{154}{3}$&$- \frac{1610}{3}$&$- \frac{34588}{5}$&$- \frac{4666292}{45}$&$- \frac{553625626}{315}$&$- \frac{1158735422}{35}$\\
\hline
\end{tabular}
}
\caption{First coefficients of the asymptotic expansions of $C_n$ and $M_n$.}
\label{tab:coefficientsCM}
\end{table}

The probability of a random chord diagram with $n$ chords to be connected is therefore $e^{-1}(1- \frac{5}{4n}) + \bigO(\frac{1}{n^2})$.

\subsection{Monolithic chord diagrams}
A chord diagram is called monolithic if it consists only of a connected component and of isolated chords which do not `contain' each other \cite{Flajolet2000}. That means with $(a,b)$ and $(c,d)$ the labels of two chords, it is not allowed that $a < c < d< b$ nor $c < a < b < d$. Let $M(x) = \sum_{n=0} M_n x^n$ be the generating function of monolithic chord diagrams. Following \cite{Flajolet2000}, $M(x)$ fulfills 
\begin{align} M(x) = C\left( \frac{x}{(1-x)^2} \right). \end{align}
Clearly, Theorem~\ref{thm:chainrule} implies that $M \in \fring{x}{2}{{\frac12}}$, because $C \in \fring{x}{2}{{\frac12}}$ and $\frac{x}{(1-x)^2} \in \R\{x\} \subset \fring{x}{2}{{\frac12}}$.
Using the $\asyOp^2_{\frac12}$-derivation on both sides of this equation together with the result for $\left(\asyOp^2_{\frac12} C\right)(x)$ in eq.\ \eqref{eqn:asyC} gives
\begin{align} \label{eqn:asyM} \begin{split} \left(\asyOp^2_{\frac12} M\right)(x) &= \frac{1}{\sqrt{2\pi}}\frac{1}{(1-x)} \frac{x}{M(x)}e^{1 - \frac{x}{2}-\frac{(1-x)^2}{2x} (2M(x) + M(x)^2)} \\ &= \frac{1}{\sqrt{2\pi}}\left(1 - 4 x -6 x^2 - \frac{154}{3} x^3 - \frac{1610}{3} x^4 - \frac{34588}{5} x^5 + \cdots \right). \end{split}  \end{align}
Some additional coefficients are given in Table \ref{tab:coefficientsCM}.
The probability of a random chord diagram with $n$ chords to be non-monolithic is therefore $1-\left(1-\frac{4}{2n-1} + \bigO(\frac{1}{n^2})\right) = \frac{2}{n} + \bigO(\frac{1}{n^2})$.

\subsection{Simple permutations}

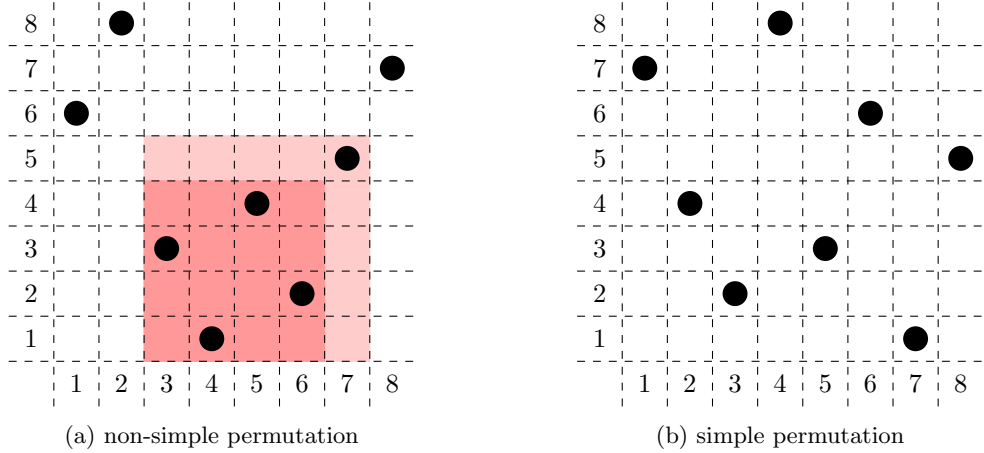
\begin{figure}%
\begin{subfigure}[b]{0.5\textwidth}%
\centering
\begin{tikzpicture}[scale=0.6] \draw[color=red!20!white,fill] (1.500000,-0.500000) rectangle (6.500000,4.500000); \draw[color=red!40!white,fill] (1.500000,-0.500000) rectangle (5.500000,3.500000); \node at (-1, 0){$1$}; \node at (0, -1){$1$}; \node at (-1, 1){$2$}; \node at (1, -1){$2$}; \node at (-1, 2){$3$}; \node at (2, -1){$3$}; \node at (-1, 3){$4$}; \node at (3, -1){$4$}; \node at (-1, 4){$5$}; \node at (4, -1){$5$}; \node at (-1, 5){$6$}; \node at (5, -1){$6$}; \node at (-1, 6){$7$}; \node at (6, -1){$7$}; \node at (-1, 7){$8$}; \node at (7, -1){$8$}; \draw[dashed] (-0.500000,-1.5)--(-0.500000,7.500000); \draw[dashed] (-1.5,-0.500000)--(7.500000,-0.500000); \draw[dashed] (0.500000,-1.5)--(0.500000,7.500000); \draw[dashed] (-1.5,0.500000)--(7.500000,0.500000); \draw[dashed] (1.500000,-1.5)--(1.500000,7.500000); \draw[dashed] (-1.5,1.500000)--(7.500000,1.500000); \draw[dashed] (2.500000,-1.5)--(2.500000,7.500000); \draw[dashed] (-1.5,2.500000)--(7.500000,2.500000); \draw[dashed] (3.500000,-1.5)--(3.500000,7.500000); \draw[dashed] (-1.5,3.500000)--(7.500000,3.500000); \draw[dashed] (4.500000,-1.5)--(4.500000,7.500000); \draw[dashed] (-1.5,4.500000)--(7.500000,4.500000); \draw[dashed] (5.500000,-1.5)--(5.500000,7.500000); \draw[dashed] (-1.5,5.500000)--(7.500000,5.500000); \draw[dashed] (6.500000,-1.5)--(6.500000,7.500000); \draw[dashed] (-1.5,6.500000)--(7.500000,6.500000); \node[circle,fill] at (0,5){}; \node[circle,fill] at (1,7){}; \node[circle,fill] at (2,2){}; \node[circle,fill] at (3,0){}; \node[circle,fill] at (4,3){}; \node[circle,fill] at (5,1){}; \node[circle,fill] at (6,4){}; \node[circle,fill] at (7,6){}; \end{tikzpicture}
\subcaption{non-simple permutation}
\end{subfigure}%
\begin{subfigure}[b]{0.5\textwidth}%
\centering
\begin{tikzpicture}[scale=0.6] \node at (-1, 0){$1$}; \node at (0, -1){$1$}; \node at (-1, 1){$2$}; \node at (1, -1){$2$}; \node at (-1, 2){$3$}; \node at (2, -1){$3$}; \node at (-1, 3){$4$}; \node at (3, -1){$4$}; \node at (-1, 4){$5$}; \node at (4, -1){$5$}; \node at (-1, 5){$6$}; \node at (5, -1){$6$}; \node at (-1, 6){$7$}; \node at (6, -1){$7$}; \node at (-1, 7){$8$}; \node at (7, -1){$8$}; \draw[dashed] (-0.500000,-1.5)--(-0.500000,7.500000); \draw[dashed] (-1.5,-0.500000)--(7.500000,-0.500000); \draw[dashed] (0.500000,-1.5)--(0.500000,7.500000); \draw[dashed] (-1.5,0.500000)--(7.500000,0.500000); \draw[dashed] (1.500000,-1.5)--(1.500000,7.500000); \draw[dashed] (-1.5,1.500000)--(7.500000,1.500000); \draw[dashed] (2.500000,-1.5)--(2.500000,7.500000); \draw[dashed] (-1.5,2.500000)--(7.500000,2.500000); \draw[dashed] (3.500000,-1.5)--(3.500000,7.500000); \draw[dashed] (-1.5,3.500000)--(7.500000,3.500000); \draw[dashed] (4.500000,-1.5)--(4.500000,7.500000); \draw[dashed] (-1.5,4.500000)--(7.500000,4.500000); \draw[dashed] (5.500000,-1.5)--(5.500000,7.500000); \draw[dashed] (-1.5,5.500000)--(7.500000,5.500000); \draw[dashed] (6.500000,-1.5)--(6.500000,7.500000); \draw[dashed] (-1.5,6.500000)--(7.500000,6.500000); \node[circle,fill] at (0,6){}; \node[circle,fill] at (1,3){}; \node[circle,fill] at (2,1){}; \node[circle,fill] at (3,7){}; \node[circle,fill] at (4,2){}; \node[circle,fill] at (5,5){}; \node[circle,fill] at (6,0){}; \node[circle,fill] at (7,4){}; \end{tikzpicture}
\subcaption{simple permutation}
\end{subfigure}%
\caption{Illustrations of simple and non-simple permutations.  The (non-trivial) intervals that map to intervals are indicated by red squares.}%
\label{fig:simple_non_simple_perm}
\end{figure}

A permutation is called simple if it does not map a non-trivial interval to another interval. Expressed formally, the permutation $\pi \in S_n^\text{simple} \subset S_n$ if and only if $\pi([i,j]) \neq [k,l]$ for all $i,j,k,l \in [1,n]$ with $2\leq |[i,j]| \leq n-1$. The difference between simple and non-simple permutations is illustrated in Figure~\ref{fig:simple_non_simple_perm}. See Albert, Atkinson and Klazar \cite{albert2003enumeration} for a detailed exposition of simple permutations. Set $S(x) = \sum_{n=4}^\infty |S_n^\text{simple}| x^n$, the generating function of \textit{simple} permutations\footnote{We adopt the convention of Albert, Atkinson and Klazar and do not consider permutations below order $4$ as simple.}, and $F(x) = \sum_{n=1}^\infty n! x^n$, the generating function of \textit{all} permutations. Following \cite{albert2003enumeration}, $S(x)$ and $F(x)$ are related by the equation
\begin{align} \label{eqn:func_S} \frac{F(x) - F(x)^2}{1+F(x)} = x + S(F(x)), \end{align}
which can be solved iteratively for the coefficients of $S(x)$:
\begin{align} S(x) &= 2 x^4 + 6 x^5 + 46 x^6 + 338 x^7 + 2926 x^8 + \cdots \end{align}
This sequence is entry \texttt{A111111} of the OEIS \cite{oeis} with the different convention, $\texttt{A111111}=x+2x^2+S(x)$.

As $n! = \G{n}{1}{1}$, $F(x) \in \fring{x}{1}{1}$ and $\left(\asyOp^1_1 F\right) = 1$ by Definitions~\ref{def:Fpowerseries} and \ref{def:basic_asymp_definition}. Therefore, the full asymptotic expansion of $S(x)$ can be obtained by applying the general chain rule to both sides of eq.\ \eqref{eqn:func_S}. Alternatively, eq.\ \eqref{eqn:func_S} implies $\frac{x-x^2}{1+x} = F^{-1}(x) + S(x)$ with $F^{-1}(F(x)) = x$. 
By Theorem~\ref{thm:chainrule}, it follows from $F \in \fring{x}{1}{1}$, $F_0 =0$ and $F_1 = 1$ that $F^{-1} \in \fring{x}{1}{1}$. By linearity and $\frac{x-x^2}{1+x} \in \R\{x\} \subset \fring{x}{1}{1}$, we also have $S \in \fring{x}{1}{1}$.
The expression for the asymptotic expansion of $F^{-1}(x)$ in terms of $\left(\asyOp^1_1 F\right)(x)$ from eq.\ \eqref{eqn:inverse_asymp} gives
\begin{align} \left(\asyOp^1_1 S\right)(x) &= \left(\asyOp^1_1 \frac{x-x^2}{1+x}\right)(x)-\left(\asyOp^1_1 F^{-1}\right)(x) = (F^{-1})'(x) \frac{x}{F^{-1}(x)} e^{\frac{1}{x}-\frac{1}{F^{-1}(x)}}, \end{align}
where $\frac{x-x^2}{1+x} \in \ker \asyOp^1_1$ was used. Observe that $F(x)$ fulfills the differential equation $x^2 F'(x) +(x-1)F(x)+x=0$, from which a non-linear differential equation for $F^{-1}(x)$ can be deduced, because $F'(F^{-1}(x)) (F^{-1})'(x) = 1$:
\begin{gather*} (F^{-1})'(x) = \frac{1}{F'(F^{-1}(x))}= \frac{F^{-1}(x)^2}{(1-F^{-1}(x))x-F^{-1}(x)}. \end{gather*}
Using this together with $\frac{x-x^2}{1+x} = F^{-1}(x) + S(x)$ gives
\begin{align} \label{eqn:asymp_genfun_Simple} \left(\asyOp^1_1 S\right)(x) &= \frac{xF^{-1}(x)}{x-(1+x)F^{-1}(x)} e^{\frac{1}{x}-\frac{1}{F^{-1}(x)}} = \frac{1}{1+x}\frac{1-x - (1+x)\frac{S(x)}{x} }{1 + (1+x) \frac{S(x)}{x^2} } e^{-\frac{2 + (1+x)\frac{S(x)}{x^2} }{1-x - (1+x)\frac{S(x)}{x}}}. \end{align}
The coefficients of $\left(\asyOp^1_1 S\right)(x)$ can be computed iteratively. The first coefficients are
\begin{align} \left(\asyOp^1_1 S\right)(x) &= e^{-2} \left( 1 - 4 x + 2 x^2 - \frac{40}{3} x^3 - \frac{182}{3} x^4 - \frac{7624}{15} x^5  + \cdots \right). \end{align}
By Definitions~\ref{def:Fpowerseries} and \ref{def:basic_asymp_definition}, this is an expression for the asymptotics of the number of simple permutations 
\begin{align} |S_n^\text{simple}| &= \sum_{k=0}^{R-1} (n-k)! [x^k] \left(\asyOp^1_1 S\right)(x) + \bigO\left( (n-R)! \right)&& \forall R \in \N_0. \end{align}
Therefore, the asymptotic expansion starts with
\begin{align*} |S_n^\text{simple}| &= e^{-2} \left( n! - 4 (n-1)! + 2 (n-2)! - \frac{40}{3} (n-3)!   - \frac{182}{3} (n-4)!  + \cdots \right).     \end{align*}
Albert, Atkinson and Klazar \cite{albert2003enumeration} calculated the first three terms of this expansion. With the presented methods the calculation of the asymptotic expansions $\left(\asyOp^1_1 S\right)(x)$ or $\left(\asyOp^1_1 F^{-1}\right)(x)$ up to order $n$ is as easy as calculating the expansion of $S(x)$ or $F^{-1}(x)$ up to order $n+2$. 
Some additional coefficients are given in Table \ref{tab:coefficientsS}.

\begin{table}
\centering
\tiny{
\def\arraystretch{1.5}
\begin{tabular}{|c||c||c|c|c|c|c|c|c|c|c|}
\hline
sequence&$0$&$1$&$2$&$3$&$4$&$5$&$6$&$7$&$8$&$9$\\
\hline\hline
$e^{2}(\asyOp^1_{1} S)$&$1$&$-4$&$2$&$- \frac{40}{3}$&$- \frac{182}{3}$&$- \frac{7624}{15}$&$- \frac{202652}{45}$&$- \frac{14115088}{315}$&$- \frac{30800534}{63}$&$- \frac{16435427656}{2835}$\\
\hline
\end{tabular}
}
\caption{First coefficients of the asymptotic expansion of $|S_n^\text{simple}|$.}
\label{tab:coefficientsS}
\end{table}

\begin{remark}
The examples above are chosen to demonstrate that given a (functional) equation which relates two power series in $\fring{x}{\alpha}{\beta}$, it is often an easy task to calculate the full asymptotic expansion of one of the power series from the asymptotic expansion of the other power series. 
Applications include functional equations for `irreducible combinatorial objects'. The two examples fall into this category. Irreducible combinatorial objects were studied in general by Beissinger \cite{beissinger1985}. 
\end{remark}
\begin{remark}
In quantum field theory  the \textit{coupling}, an expansion parameter, needs to be re\-para\-met\-rized in the process of \textit{renormalization} \cite{connes2001renormalization}. Those reparametrizations are merely compositions of power series which are believed to be \textit{Gevrey-1}. Theorem~\ref{thm:chainrule} might be useful for the resummation of renormalized quantities in quantum field theory.
Dyson-Schwinger equations in quantum field theory can be stated as functional equations of a form similar to the above \cite{Broadhurst2001,borinsky2015algebraic}. These considerations were the subject of the publication \cite{borinsky2017}, where the presented formalism was applied to zero-dimensional quantum field theory and the enumeration of graphs.
\end{remark}
\begin{remark}
Eqs. \eqref{eqn:asyC}, \eqref{eqn:asyM} and \eqref{eqn:asymp_genfun_Simple} expose another interesting algebraic property. Proposition~\ref{prop:betashiftlow} and the chain rule imply that $(\asyOp^2_{\frac12} C)(x)\in\fring{x}{2}{\frac32}$, $(\asyOp^2_{\frac12} M)(x) \in \fring{x}{2}{\frac32}$ and $(\asyOp^1_1 S)(x) \in \fring{x}{1}{3}$. This way, the `higher-order' asymptotics of the asymptotic sequence can be calculated by iterating the application of the $\asyOp$ map. With resurgence, it might be possible to construct \textit{convergent} large-order expansions for these cases.
The fact that the asymptotics of each sequence may be expressed as a combination of polynomial and exponential expressions of the original sequence can be seen as an avatar of resurgence. 
\end{remark}

\section*{Acknowledgements}
Many thanks to Dirk Kreimer for steady encouragement and counseling. I wish to express my gratitude to David Broadhurst. He sparked my interest in asymptotic expansions and I benefited greatly from our discussions. I wish to thank In\^es Aniceto who patiently introduced me to the basics of resurgence theory, David Sauzin for encouragement and patiently clarifying details of the connections to resurgence, Julien Courtiel for helpful comments regarding the closure properties of $\fring{x}{\alpha}{\beta}$ under functional inversion and Karen Yeats for hospitality during the revision of this article. I also wish to thank the anonymous referee for carefully reading the manuscript and providing many helpful comments.

\bibliographystyle{plain}
\bibliography{literature.bib}

\begin{thebibliography}{10}

\bibitem{albert2003enumeration}
MH~Albert, MD~Atkinson, and M~Klazar.
\newblock \titleurl{The enumeration of simple
  permutations}{http://emis.ams.org/journals/JIS/VOL6/Albert/albert.pdf}.
\newblock {\em Journal of Integer Sequences}, volume 6:03.4.4, 2003.

\bibitem{alvarez2004langer}
G~Alvarez.
\newblock \titledoi{{Langer--Cherry} derivation of the multi--instanton
  expansion for the symmetric double well}{10.1063/1.1767988}.
\newblock {\em Journal of mathematical physics}, 45(8):3095--3108, 2004.

\bibitem{aniceto2011resurgence}
I~Aniceto, R~Schiappa, and M~Vonk.
\newblock \titledoi{The resurgence of instantons in string
  theory}{10.4310/CNTP.2012.v6.n2.a3}.
\newblock {\em Communications in Number Theory and Physics}, 6(2):339--496,
  2012.

\bibitem{beissinger1985}
JS~Beissinger.
\newblock \titledoi{The enumeration of irreducible combinatorial
  objects}{10.1016/0097-3165(85)90065-2}.
\newblock {\em Journal of Combinatorial Theory, Series A}, 38(2):143--169,
  1985.

\bibitem{bender1969anharmonic}
CM~Bender and TT~Wu.
\newblock \titledoi{Anharmonic oscillator}{10.1103/PhysRev.184.1231}.
\newblock {\em Physical Review}, 184(5):1231--1260, 1969.

\bibitem{bender1974asymptotic}
EA~Bender.
\newblock \titledoi{Asymptotic methods in enumeration}{10.1137/1016082}.
\newblock {\em SIAM Review}, 16(4):485--515, 1974.

\bibitem{bender1975asymptotic}
EA~Bender.
\newblock \titledoi{An asymptotic expansion for the coefficients of some formal
  power series}{10.1112/jlms/s2-9.3.451}.
\newblock {\em Journal of the London Mathematical Society}, 2(3):451--458,
  1975.

\bibitem{bender1978asymptotic}
EA~Bender and ER~Canfield.
\newblock \titledoi{The asymptotic number of labeled graphs with given degree
  sequences}{10.1016/0097-3165(78)90059-6}.
\newblock {\em Journal of Combinatorial Theory, Series A}, 24(3):296--307,
  1978.

\bibitem{Bender1984}
EA~Bender and LB~Richmond.
\newblock \titledoi{An asymptotic expansion for the coefficients of some power
  series {II}: {Lagrange} inversion}{10.1016/0012-365X(84)90043-8}.
\newblock {\em Discrete Mathematics}, 50:135--141, 1984.

\bibitem{bergeron1998combinatorial}
F~Bergeron, G~Labelle, and P~Leroux.
\newblock {\em Combinatorial species and tree-like structures}, volume~67.
\newblock Cambridge University Press, 1998.

\bibitem{borinsky2015algebraic}
M~Borinsky.
\newblock \titledoi{Algebraic lattices in qft
  renormalization}{10.1007/s11005-016-0843-9}.
\newblock {\em Letters in Mathematical Physics}, 106(7):879--911, 2016.

\bibitem{borinsky2017}
M~Borinsky.
\newblock \titledoi{Renormalized asymptotic enumeration of {Feynman}
  diagrams}{10.1016/j.aop.2017.07.009}.
\newblock {\em Annals of Physics}, 385:95--135, 2017.

\bibitem{borinsky2016generating}
M~Borinsky.
\newblock \titleurl{Generating asymptotics for factorially divergent sequences
  (extended
  abstract)}{http://www.mat.univie.ac.at/~slc/wpapers/FPSAC2017/28\%20Borinsky.html}.
\newblock In {\em S\'eminaire Lotharingien de Combinatoire -- Proceedings to
  FPSAC 2017}, volume 78B, page 12pp., 2017.

\bibitem{Broadhurst2001}
DJ~Broadhurst and D~Kreimer.
\newblock \titledoi{Exact solutions of {Dyson--Schwinger} equations for
  iterated one-loop integrals and propagator-coupling
  duality}{10.1016/S0550-3213(01)00071-2}.
\newblock {\em Nuclear Physics B}, 600(2):403--422, 2001.

\bibitem{connes2001renormalization}
A~Connes and D~Kreimer.
\newblock \titledoi{Renormalization in quantum field theory and the
  {Riemann--Hilbert} problem {II}: The $\beta$-function, diffeomorphisms and
  the renormalization group}{10.1007/s002200050779}.
\newblock {\em Communications in Mathematical Physics}, 216(1):215--241, 2001.

\bibitem{bruijn1970asymptotic}
NG~de~Bruijn.
\newblock {\em Asymptotic Methods in Analysis}.
\newblock Bibliotheca mathematica. Dover Publications, 1970.

\bibitem{dunne2012resurgence}
GV~Dunne and M~{\"U}nsal.
\newblock \titledoi{Resurgence and trans-series in quantum field theory: the
  $\mathbb{C}\mathbb{P}^{N-1}$ model}{10.1007/JHEP11(2012)170}.
\newblock {\em Journal of High Energy Physics}, 2012(11):1--86, 2012.

\bibitem{ecalle1981fonctions}
J~{\'E}calle.
\newblock Les fonctions r{\'e}surgentes.
\newblock {\em Publ. math. d'Orsay/Univ. de Paris, Dep. de math.}, 1981.

\bibitem{Flajolet2000}
P~Flajolet and M~Noy.
\newblock \titledoi{Analytic combinatorics of chord
  diagrams}{10.1007/978-3-662-04166-6_17}.
\newblock In {\em Formal Power Series and Algebraic Combinatorics}, pages
  191--201. Springer, 2000.

\bibitem{flajolet2009analytic}
P~Flajolet and R~Sedgewick.
\newblock {\em Analytic combinatorics}.
\newblock Cambridge University Press, 2009.

\bibitem{garoufalidis2012asymptotics}
S~Garoufalidis, A~Its, A~Kapaev, and M~Marino.
\newblock \titledoi{Asymptotics of the instantons of {Painlev}{\'e}
  {I}}{10.1093/imrn/rnr029}.
\newblock {\em International Mathematics Research Notices}, 2012(3):561--606,
  2012.

\bibitem{hsieh2012basic}
PF~Hsieh and Y~Sibuya.
\newblock {\em Basic Theory of Ordinary Differential Equations}.
\newblock Universitext. Springer New York, 2012.

\bibitem{kleitman1970proportions}
DJ~Kleitman.
\newblock \titledoi{Proportions of irreducible
  diagrams}{10.1002/sapm1970493297}.
\newblock {\em Studies in Applied Mathematics}, 49(3):297--299, 1970.

\bibitem{le2012large}
J~Le~Guillou and J~Zinn-Justin.
\newblock {\em Large-order behaviour of perturbation theory}, volume~7.
\newblock Elsevier, 2012.

\bibitem{mitschi2016divergent}
C~Mitschi and D~Sauzin.
\newblock {\em Divergent Series, Summability and Resurgence I}.
\newblock Springer, 2016.

\bibitem{oeis}
NJA Sloane.
\newblock The on-line encyclopedia of integer sequences. http://oeis.org/,
  2005.

\bibitem{stein1978class}
PR~Stein and CJ~Everett.
\newblock \titledoi{On a class of linked diagrams {II.}
  asymptotics}{10.1016/0012-365X(78)90162-0}.
\newblock {\em Discrete Mathematics}, 21(3):309--318, 1978.

\bibitem{wilf2013generatingfunctionology}
HS~Wilf.
\newblock {\em generatingfunctionology}.
\newblock Elsevier, 2013.

\end{thebibliography}

\end{document}